\documentclass[11pt,a4paper]{article}
\usepackage[left=25mm,right=25mm,top=24mm,bottom=24mm]{geometry}
\usepackage{amsmath,amssymb,mathtools,booktabs,array,needspace,microtype}
\usepackage[hidelinks,unicode,hypertexnames=false,hyperfootnotes=false]{hyperref}
\usepackage{bookmark}
\allowdisplaybreaks[2]
\newcommand{\HH}{\mathbb H}
\newcommand{\RR}{\mathbb R}
\newcommand{\dd}{\,\mathrm d}
\newcommand{\norm}[1]{\left\lVert#1\right\rVert}

\newcommand{\result}[1]{\par\addvspace{9pt}\Needspace{5\baselineskip}\noindent\textbf{#1}\ }
\newcommand{\step}[1]{\par\addvspace{5pt}\Needspace{3\baselineskip}\noindent\textit{#1}\ }
\newcommand{\qedbox}{\hfill$\square$}

\newcommand{\Py}{\mathcal P_y}

\newcommand{\Tk}{\mathcal T_{|k|}}
\newcommand{\Tl}{\mathcal T_\lambda}
\newcommand{\Pil}{\Pi_\ell}
\newcommand{\TT}{\mathcal T}
\newcommand{\PP}{\mathcal P}
\DeclareMathOperator{\curl}{curl}
\DeclareMathOperator{\diver}{div}
\title{A Liouville theorem for steady
Navier--Stokes flows in the upper half-plane}

\author{Yibing Sun\textsuperscript{$*$}
\and Guo Zhang\textsuperscript{$\dagger$}}
\newcommand{\firstauthoraddress}{School of Mathematical Sciences, Beijing Normal University, Haidian District Xinjiekou Waidajie
No.19, 100875, Beijing, P. R. China.
Email: \texttt{yibingsun@mail.bnu.edu.cn}.}
\newcommand{\secondauthoraddress}{Corresponding author, Hubei Key Lab-Math.Sci.,
Central China Normal University, Wuhan 430079, P. R. China.
Email: \texttt{guozhang625@hotmail.com}.}
\newcommand{\printauthoraddresses}{%
  \begingroup
  \renewcommand{\thefootnote}{\fnsymbol{footnote}}%
  \footnotetext[1]{\raggedright\firstauthoraddress}%
  \footnotetext[2]{\raggedright\secondauthoraddress}%
  \endgroup
}
\date{} % Insert a date here if desired.
\makeatletter
\renewcommand{\@maketitle}{%
  \newpage\null\vskip 8pt
  \begin{center}%
    {\fontsize{20.74}{25}\selectfont\@title\par}%
    \vskip 16pt
    {\large\lineskip .5em
      \begin{tabular}[t]{c}\@author\end{tabular}\par}%
    \ifx\@date\@empty\else\vskip 10pt{\normalsize\@date\par}\fi
  \end{center}%
  \par\vskip 8pt}
\makeatother
\renewenvironment{abstract}{%
  \par\begin{center}\small\bfseries\abstractname\end{center}%
  \begin{list}{}{\setlength{\leftmargin}{24pt}\setlength{\rightmargin}{24pt}%
    \setlength{\topsep}{3pt}\setlength{\parsep}{0pt}}%
  \item\relax\small\linespread{1.08}\selectfont
}{\end{list}\vspace{6pt}}
\hypersetup{pdftitle={A Liouville theorem for steady Navier-Stokes flows in the upper half-plane},pdfauthor={Guo Zhang},pdfsubject={Finite-Dirichlet-integral solutions in the upper half-plane},pdfkeywords={Navier-Stokes, Liouville theorem, half-plane, boundary vorticity, Oseen operator}}
\AtBeginDocument{\synctex=-1\relax}
\begin{document}
% sections123
\maketitle
\printauthoraddresses
\begin{abstract}
We study steady incompressible Navier--Stokes flows in the upper half-plane with nonzero constant boundary velocity. We prove that every solution with finite Dirichlet integral is the corresponding constant flow. The proof distinguishes outflow, inflow, and tangential boundary velocities. In the outflow case, a vorticity identity yields a weighted estimate and then square integrability of the perturbation. In the inflow case, weighted pressure estimates allow the nonlinear terms to be absorbed in an energy identity. For tangential flows, the momentum equations and the boundary conditions on the stream function give complementary bulk and boundary-vorticity identities. Uniform weighted estimates for the associated Riesz transforms and Oseen response yield square integrability of the normal velocity component. A truncated stream-function test then proves rigidity.
\end{abstract}

\section{Introduction}
We study the Liouville problem for the steady incompressible Navier--Stokes equations in the upper half-plane
\[
\HH=\{(x,y)\in\RR^2:y>0\}.
\]
More precisely, we consider
\begin{equation}
\left\{
\begin{aligned}
&-\Delta U+(U\cdot\nabla)U+\nabla P=0,\ \ \ \diver U=0 \ \ \ \text{in }\HH,\\
&U\big|_{\partial\HH}=B,\ \ \ \ \ \ E:=\int_\HH|\nabla U|^2\dd x\dd y<\infty,
\end{aligned}
\right.
\tag{1.1}\label{eq:intro-system}
\end{equation}
where $B=(B_1,B_2)$ is a nonzero constant vector. A solution with finite Dirichlet integral is usually called a $D$-solution. We ask whether every $D$-solution of \eqref{eq:intro-system} is the constant flow $U\equiv B$, without prescribing its behavior at infinity.

In the whole plane, Gilbarg and Weinberger \cite{GW1978} proved that every $D$-solution is constant. The work of Koch, Nadirashvili, Seregin, and \v Sver\'ak \cite{KNSS2009} on bounded ancient solutions also implies that every bounded steady solution in $\RR^2$ is constant. Fuchs and Zhong \cite{FZ2011} extended this conclusion to solutions satisfying $|U(z)|=O(|z|^\alpha)$ for some $0\le\alpha<1/7$. Bildhauer, Fuchs, and Zhang \cite{BFZ2013} subsequently enlarged this range to $0\le\alpha<1/3$. Under an integrability assumption, Zhang \cite{Zhang2015} proved that an entire steady solution with $U\in L^p(\RR^2)$ for some $1<p<\infty$ must vanish. Related Liouville theorems for non-Newtonian fluids were obtained by Fuchs \cite{Fuchs2012}, Fuchs and Zhang \cite{FuchsZhang2012}, and Bildhauer, Fuchs, and Zhang \cite{BFZ2013}; see also \cite{Zhang2015} for finite-energy and integrability criteria for shear-thickening fluids.

In the half-plane, the boundary condition has a decisive role. For homogeneous Dirichlet data, Seregin \cite{Seregin2014} proved that a bounded mild ancient solution with uniformly bounded kinetic energy is zero. In particular, this result excludes stationary solutions in that class whose velocities belong to $L^2(\HH)\cap L^\infty(\HH)$. Giga, Hsu, and Maekawa \cite{GHM2014} obtained a Liouville theorem for ancient solutions with one-signed vorticity, under uniform H\"older bounds, a prescribed pressure decomposition, and a type-I temporal decay condition.  Under the perfect-slip conditions $\partial_yU_1(x,0)=0,\ U_2(x,0)=0,$
Wang and Wang \cite{WW2023,WW2022} proved constancy for half-plane $D$-solutions. In that setting the vorticity vanishes on the boundary, which allows a direct use of the maximum principle.  For the forced equation with homogeneous boundary data, Calderon, Le, and Phan \cite{CLP2024} established existence and weak--strong uniqueness results using weighted estimates and stream functions.

Liouville questions are closely related to the behavior of the velocity and vorticity at infinity. For flow past a bounded obstacle, Leray \cite{Leray1933,Leray1934} introduced the invading-domain method to construct solutions with finite Dirichlet integral. Gilbarg and Weinberger \cite{GW1978} obtained fundamental asymptotic estimates, including convergence to a constant velocity in an angular mean-square sense. Amick \cite{Amick1988} proved uniform convergence for symmetric flows. Korobkov, Pileckas, and Russo \cite{KPR2019} later established uniform convergence at infinity for every planar exterior-domain $D$-solution, without symmetry or smallness assumptions. Identifying this limiting vector with the prescribed velocity in Leray's construction is an open question; Korobkov and Ren \cite{KR2022} resolved it at sufficiently small Reynolds numbers.

When the limiting velocity $U_\infty$ is nonzero, sharper asymptotics are available. Smith \cite{Smith1965} showed that an initial estimate
\[
|U(z)-U_\infty|=O(|z|^{-1/4-\varepsilon}),\qquad \varepsilon>0,
\]
leads to the classical wake asymptotics, with an $O(|z|^{-1/2})$ leading velocity correction in the downstream wake. Sazonov \cite{Sazonov1999} showed that an exterior-domain $D$-solution converging uniformly to $U_\infty\ne0$ satisfies such an initial algebraic estimate and hence belongs to the class of physically reasonable solutions. These results concern a neighborhood of infinity outside a bounded obstacle. Their decay conclusions do not extend to arbitrary half-plane $D$-solutions without additional boundary information.

An explicit example illustrates this distinction. Let $z=x+iy$, fix $B\in\RR^2$, and define
\[
F(z)=B_1-iB_2+\frac1{\operatorname{Log}(z+ie)},\qquad
U=(\Re F,-\Im F),\qquad
P=\frac{|B|^2-|U|^2}{2},
\]
where $\operatorname{Log}$ is the principal branch of the complex logarithm. Since $\Im(z+ie)\ge e$ on $\overline\HH$, the function $F$ is holomorphic in a neighborhood of $\overline\HH$ and $\Re\operatorname{Log}(z+ie)\ge1$. The Cauchy--Riemann equations give $\diver U=\curl U=0$ and $\Delta U=0$. Thus $(U\cdot\nabla)U=\nabla(|U|^2/2)$, and $(U,P)$ is a smooth steady Navier--Stokes solution. Moreover,
\[
\int_\HH|\nabla U|^2
=2\int_\HH\frac{\dd x\dd y}{|z+ie|^2|\operatorname{Log}(z+ie)|^4}
\le2\pi\int_e^\infty\frac{\dd r}{r(\log r)^4}
=\frac{2\pi}{3}.
\]
Nevertheless,
\[
U(z)\longrightarrow B,\qquad
|U(z)-B|=\frac1{|\operatorname{Log}(z+ie)|}
\sim\frac1{\log|z|}
\quad\text{as }|z|\to\infty,
\]
uniformly in $\HH$. Hence $U-B$ has no algebraic decay rate. Its boundary trace is not constant, so this example does not satisfy \eqref{eq:intro-system}. It shows why finite Dirichlet integral and a nonzero far-field limit alone cannot supply the exterior-domain decay estimates in the half-plane.

For constant boundary data, we instead use weighted estimates to obtain the integrability needed for an energy argument. Our main result is the following.

\result{Theorem 1.1 (Liouville theorem).}
Let $B\in\RR^2\setminus\{0\}$. Suppose that
$U\in H^1_{\rm loc}(\overline\HH;\RR^2)$ and
$P\in L^2_{\rm loc}(\overline\HH)$ satisfy \eqref{eq:intro-system}. Then $U\equiv B$ in $\HH$, and $P$ is constant.

Set $W=U-B=(u,v)$. Then \eqref{eq:intro-system} is equivalent to
\begin{equation}
\left\{
\begin{aligned}
&-\Delta W+(B+W)\cdot\nabla W+\nabla P=0, \ \ \ \diver W=0\ \ \  \text{in }\HH,\\
&W\big|_{\partial\HH}=0,\ \ \ \ \  E=\int_\HH|\nabla W|^2\dd x\dd y<\infty.
\end{aligned}
\right.
\tag{1.2}\label{eq:intro-perturbation}
\end{equation}
Thus it remains to prove $W\equiv0$. We use the notation
\begin{equation}
\omega=v_x-u_y,\qquad
W=\nabla^\perp\psi:=(\psi_y,-\psi_x),\qquad
\omega=-\Delta\psi,
\tag{1.3}\label{eq:intro-vorticity}
\end{equation}
with $\psi(x,0)=\psi_y(x,0)=0$. The Hardy inequalities and the interior velocity oscillation estimate of Korobkov and Ren \cite{KR2025} first yield
\[
W/y,\ \psi/y^2,\ \omega\in L^2(\HH),\qquad
\|W\|_\infty\le C\sqrt E,\qquad
W(z)\longrightarrow0\quad\text{uniformly as }|z|\to\infty.
\]
The proof then separates into three cases according to the sign of $B_2$.

\smallskip
\noindent\textit{Outflow: $B_2<0$.}
Write $B_2=-\gamma$, with $\gamma>0$, and set
\[
g(y)=\int_\RR\omega^2\dd x,\qquad
d(y)=\int_\RR|\nabla\omega|^2\dd x,\qquad
n(y)=\int_\RR v\omega^2\dd x.
\]
A horizontal-slice calculation in the vorticity equation gives
\[
g'(y)=-\gamma g(y)+n(y)-2\int_y^\infty d(s)\dd s.
\]
Since $v$ tends uniformly to zero, this identity gives exponential decay of $g$ and, in particular, $y\omega\in L^2(\HH)$. A weighted estimate for $-\Delta\psi=\omega$ then implies
\[
\|W\|_2\le C\|y\omega\|_2.
\]
With $W\in L^2(\HH)$ established, a cutoff energy test in \eqref{eq:intro-perturbation} gives $E=0$.

\smallskip
\noindent\textit{Inflow: $B_2>0$.}
For $\rho>0$, introduce the even weights
\[
\chi_\rho(y)=\frac\rho{\rho+|y|},\qquad
\vartheta_\rho(y)=\chi_\rho(y)^2,\qquad
w_\rho(y)=\frac{\vartheta_\rho(y)}\rho
=\frac\rho{(\rho+|y|)^2}.
\]
On $y>0$, $w_\rho=-\chi_\rho'$. Put $A_\rho=\int_\HH w_\rho|W|^2$. Decompose the pressure as $P=\pi+p_s$, where $\pi$ is the quadratic Riesz-transform pressure and $p_s$ is the harmonic pressure determined by the boundary vorticity. Testing with $\chi_\rho W$ gives
\[
\begin{aligned}
\int_\HH\chi_\rho|\nabla W|^2+\frac{B_2}{2}A_\rho
={}&\frac12\int_\HH\chi_\rho''|W|^2
+\frac12\int_\HH\chi_\rho'v|W|^2\\
&-\int_\HH\nabla\pi\cdot(\chi_\rho W)
+\int_\HH\chi_\rho'p_sv.
\end{aligned}
\]
The main difficulty is the quadratic pressure term. Reflection identities reduce its estimate to weighted Riesz-transform bounds, leading to
\[
\left|\int_\HH\nabla\pi\cdot(\chi_\rho W)\right|
\le C\|W\|_\infty A_\rho.
\]
To remove any smallness assumption on $W$, we split $W^e\otimes W^e$ into a fixed compact contribution and a tail where $W$ is small. The compact contribution tends to zero as $\rho\to\infty$, whereas the tail is bounded by $C\varepsilon A_\rho$. The remaining terms are controlled in the same weighted identity. The positive normal-transport term then permits absorption and yields $E=0$.

\smallskip
\smallskip
\noindent\textit{Tangential flow: $B_2=0$.}
After a horizontal reflection, write $B=(b,0)$ with $b>0$.
This case presents a different difficulty: the positive normal-transport
term in the inflow argument disappears, while finite Dirichlet integral
alone does not provide the velocity integrability needed to justify a
global energy test. In particular, the usual Calder\'on--Zygmund estimate
for the quadratic pressure involves $\|W\|_4^2$, or
$\|W\|_\infty\|W\|_2$, neither of which is available at this stage.

Fourier analysis of the tangential Oseen operator has been used
successfully to study flows near a wall. In particular,
Boeckle and Wittwer \cite{BW2012} obtained quantitative decay estimates
for solutions driven by sufficiently small compactly supported forces,
using function spaces adapted to the Fourier representation of the flow.
Our purpose here is to obtain an a priori rigidity result for arbitrary
$D$-solutions. The argument must therefore close without assuming a
small Dirichlet integral or a quantitative decay rate at infinity.

The starting point is an exact identity coupling the normal velocity
to the interior and boundary vorticity:
\[
bv+G(W\otimes W)=\omega-h_f,
\qquad
f=\omega(\cdot,0),\qquad h_f=\Py f,
\]
where
\[
G(W\otimes W)
=R_xR_y(u^2-v^2)+(R_y^2-R_x^2)(uv).
\]
The diagonal tensor entries are extended evenly across the boundary,
and the off-diagonal entries oddly.
Although $\omega\in L^2(\HH)$ is already known, neither term on the
right-hand side of
$bv=\omega-h_f-G(W\otimes W)$ can simply be estimated by standard
velocity bounds. The bulk term contains the tangential square $u^2$,
whereas the harmonic term has a low-frequency obstruction:
\[
\|\Py f\|_{L^2(\HH)}^2
=\frac12\int_\RR\frac{|\widehat f(k)|^2}{|k|}\dd k.
\]
Thus the available boundary regularity of $f$ does not by itself
control its Poisson extension in $L^2(\HH)$.

We address both difficulties in the same weighted norm,
\[
V_\rho(v):=\|\chi_\rho v\|_2,
\qquad
\chi_\rho(y)=\frac{\rho}{\rho+|y|}.
\]
For each fixed $\rho$, the Hardy inequality ensures that
$V_\rho(v)<\infty$.
The essential point is to obtain estimates uniform in $\rho$
whose right-hand sides involve this normal component, without
requiring the corresponding unweighted norm of the full velocity.

The bulk estimate rests on two complementary features of the
Riesz transforms, developed in Appendices A and B.
First, incompressibility allows the tangential square to be
estimated through the normal component. A dyadic decomposition
in the normal frequency, based on the classical
Littlewood--Paley and Bony framework
\cite{Bony1981,BCD2011}, yields
\[
\|R_xR_y(r^2)\|_{L^2(\RR^2)}
\le C\|\nabla r\|_{L^2(\RR^2)}
      \|s\|_{L^2(\RR^2)},
\qquad r_x+s_y=0,
\]
for divergence-free $H^1$ fields.
The derivative transfer supplied by incompressibility is crucial:
it replaces the unavailable integrability of the tangential component
by that of the normal component.
Second, the density $\vartheta_\rho=\chi_\rho^2$ lies outside the
Muckenhoupt $A_2$ class, so the classical weighted
Calder\'on--Zygmund theory \cite{Grafakos2014} does not apply directly.
We instead exploit the parity imposed by reflection.
For sources well above the observation point, cancellation between
the reflected kernels provides an additional ratio of heights,
making the weighted interactions between distant layers summable.
Together with divergence-preserving localization of the stream
function and an estimate for odd harmonic functions, this gives
\[
\|\chi_\rho G(A\otimes A)\|_2
\le C\bigl(\|\nabla A\|_2+\|A\|_\infty\bigr)V_\rho(s),
\qquad A=(r,s).
\]

The boundary contribution requires a separate use of the no-slip
condition. Both conditions $\psi=\psi_y=0$ on $y=0$ enter the
Fourier representation of the vorticity equation and determine
\[
f=T_b(W\otimes W).
\]
Consequently, the boundary vorticity is constrained by the interior
quadratic tensor, rather than being arbitrary boundary data.
At low horizontal frequencies, the response involves two distinct
normal scales, governed by $|k|$ and
$\lambda(k)=\sqrt{|k|^2+ibk}$, with $\Re\lambda>0$ and
$|\lambda(k)|\asymp\sqrt{b|k|}$.
Controlling their interaction in the finite-Dirichlet-integral
class is the main issue addressed by Lemma 5.5.

More precisely, if $\Pi_\ell$ denotes the projection onto
$0<|k|\le b/4$, Lemma 5.5 proves, for $\rho\ge b^{-1}$,
\[
\begin{aligned}
\|\chi_\rho\Py\Pi_\ell T_b(A\otimes A)\|_2
\le{}&
Cb^{-1}\|\nabla A\|_2^2\\
&+C\bigl(\|\nabla A\|_2+\|A\|_\infty\bigr)V_\rho(s)
+\frac{C\|\nabla A\|_2}{\sqrt{b\rho}}V_\rho(s).
\end{aligned}
\]
Its proof preserves the algebraic cancellations in the response
formula before estimating individual terms.
Using incompressibility, the apparently uncontrolled $r^2$
contributions are rewritten in terms of $r_y s$, $rs$, and $s^2$.
A cutoff of the stream function then separates a field supported
at finite height, to which Lemma 5.4 applies, from sources on
dyadic strips above height $\rho$.
The Fourier kernel estimates on those strips retain both the
Poisson and Oseen scales. Their summation produces the factor
$(b\rho)^{-1/2}$ and, crucially, the same quantity $V_\rho(s)$
as in the bulk estimate.
This extends the boundary estimate beyond the class in which
the normal velocity is already square integrable, thereby
avoiding a circular argument.

Finally, the uniformity of these two estimates allows us to remove
a global smallness restriction.
We decompose $W=W_c+W_\infty$ through its stream function, with
$W_c$ compactly supported and $W_\infty$ small in both
$L^\infty$ and the Dirichlet seminorm.
The compact contribution is fixed before $\rho$ varies, while
the tail contributions can be absorbed in the momentum identity.
This yields
\[
bV_\rho(v)
\le \|\omega\|_2+C_R+C\varepsilon V_\rho(v),
\qquad \rho\ge b^{-1}.
\]
After absorption, letting $\rho\to\infty$ gives
$v\in L^2(\HH)$.
This single-component integrability is sufficient:
a truncated stream-function test in the vorticity equation,
using $W\cdot\nabla\psi=0$, then gives $E=0$.
The tangential argument thus combines a normal-component estimate
for the nonlinear bulk term with a compatible low-frequency
boundary estimate, obtaining rigidity without first deriving
an algebraic far-field expansion.

Section 2 establishes the basic estimates and the pressure decomposition. Sections 3 and 4 treat outflow and inflow. Sections 5 and 6 develop the tangential argument and complete the proof of Theorem 1.1. Appendices A and B contain the anisotropic and weighted Riesz-transform estimates used in the proof.

\subsection*{Notation}
An unqualified $L^p$ norm is taken over the natural domain of the function in question; other domains are indicated explicitly. We use the unitary Fourier transform in the horizontal variable, so that $\widehat{\partial_x f}(k)=ik\widehat f(k)$. The symbols $R_x,R_y$ denote the whole-plane Riesz transforms, and $H_x$ denotes the horizontal Hilbert transform with symbol $-i\operatorname{sgn}k$. The Poisson extension is defined by $\widehat{\Py f}(k)=e^{-|k|y}\widehat f(k)$. The letter $C$ denotes a positive constant that may change from line to line. Constants denoted by $C_b$ may depend on the fixed tangential velocity $b$, but not on weight parameters or cutoff heights.

\section{Basic estimates and pressure decomposition}
Throughout this section, $B$ is an arbitrary constant vector,
all unqualified norms are taken over $\HH$, and
$E=\int_\HH |\nabla W|^2$.
Some of the results below are well known; we include their statements
and brief proofs for completeness.

\result{Lemma 2.1 }
Suppose that $\diver W=0$, $W(x,0)=0$, and $E<\infty$. Then

\[
\psi(x,y):=\int_0^yu(x,t)\dd t
\]

 satisfies $\psi_y=u$, $\psi_x=-v$, $-\Delta\psi=\omega$, and $\psi(x,0)=\psi_y(x,0)=0$. Moreover,

\begin{equation}
\norm{W/y}_2\le2\norm{\partial_yW}_2,\qquad
\norm{\psi/y^2}_2\le\frac23\norm{u/y}_2,
\tag{2.1a}
\end{equation}

\begin{equation}
\norm{\omega}_2^2=\int_{\HH}|\nabla W|^2=E.
\tag{2.1b}
\end{equation}

\step{Proof. Step 1: the second Hardy inequality.}
 Since $u(x,0)=0$ and $\psi_y=u$, for almost every $x$, integration by parts and the Cauchy--Schwarz inequality give

\[
\psi(x,\varepsilon)=\int_0^\varepsilon(\varepsilon-t)u_y(x,t)\dd t,
\qquad
\frac{|\psi(x,\varepsilon)|^2}{\varepsilon^3}
\le\frac13\int_0^\varepsilon|u_y(x,t)|^2\dd t\longrightarrow0.
\]

Integrating once more, now over $(\varepsilon,R)$, we obtain

\[
\begin{aligned}
\int_\varepsilon^R\frac{\psi^2}{y^4}\dd y
&=-\frac{\psi(x,R)^2}{3R^3}
 +\frac{\psi(x,\varepsilon)^2}{3\varepsilon^3}
 +\frac23\int_\varepsilon^R\frac{\psi u}{y^3}\dd y\\
&\le\frac{\psi(x,\varepsilon)^2}{3\varepsilon^3}
 +\frac23\left(\int_\varepsilon^R\frac{\psi^2}{y^4}\dd y\right)^{1/2}
 \left(\int_\varepsilon^R\frac{u^2}{y^2}\dd y\right)^{1/2}.
\end{aligned}
\]

Letting $\varepsilon\downarrow0$, we obtain

\[
\int_0^R\frac{\psi^2}{y^4}\dd y
\le\frac23\left(\int_0^R\frac{\psi^2}{y^4}\dd y\right)^{1/2}
\left(\int_0^R\frac{u^2}{y^2}\dd y\right)^{1/2}
.
\]

Dividing by the square root of the left-hand side, when it is nonzero, and then squaring, letting $R\to\infty$, and integrating in $x$, we obtain $\|\psi/y^2\|_2\le(2/3)\|u/y\|_2$.

\step{Step 2: the vorticity energy identity.}
For a smooth compactly supported vector field $F$, the standard identity reads

\[
\|\nabla F\|_{L^2(\RR^2)}^2
=\|\curl F\|_{L^2(\RR^2)}^2+\|\diver F\|_{L^2(\RR^2)}^2.
\]

By density, it also holds for compactly supported $H^1$ vector fields. Extend $W$ by zero to the lower half-plane and denote the extension by $\widetilde W$. Its zero trace implies

\[
\diver\widetilde W=0,\qquad \curl\widetilde W=\widetilde\omega,
\qquad \|\nabla\widetilde W\|_{L^2(\RR^2)}^2=E.
\]

Choose a radial cutoff $\eta_R$ equal to one on $B_R$ and zero outside $B_{2R}$, with $|\nabla\eta_R|\le C/R$. The first Hardy inequality gives

\[
\|\widetilde W\otimes\nabla\eta_R\|_2^2
\le\frac C{R^2}\int_{(B_{2R}\setminus B_R)\cap\HH}|W|^2
\le C\int_{\{|z|>R\}\cap\HH}|W/y|^2\longrightarrow0.
\]

Consequently, $F_R=\eta_R\widetilde W\in H_c^1(\RR^2)$ satisfies

\[
\nabla F_R\to\nabla\widetilde W,\qquad
\diver F_R\to0,\qquad \curl F_R\to\widetilde\omega
\quad\text{in }L^2.
\]

Applying the preceding identity to $F_R$ and passing to the limit proves (2.1b).\qedbox

\result{Lemma 2.2 }
 $W\in C_0(\overline\HH)$, and

\begin{equation}
\|W\|_{C_0(\overline\HH)}\le C\|\nabla W\|_{L^2(\HH)}.
\tag{2.2}
\end{equation}

\step{Proof.}
Continuity follows from standard local regularity for the planar steady equations. For $z=(x,h)\in\HH$, set

\[
D_z=B_{h/2}(z),\qquad D'_z=B_{h/4}(z),\qquad
m_z=|D'_z|^{-1}\int_{D'_z}W(w)\dd w.
\]

Apply \cite{KR2025} to the original velocity $U=B+W$, with $\delta=1/2$, and average over $D'_z$. This gives

\begin{equation}
|W(z)-m_z|\le C\norm{\nabla W}_{L^2(D_z)}.
\tag{2.3}
\end{equation}

with a constant independent of $z$, $h$, and $B$. On $D'_z$ we have $3h/4<y<5h/4$ and $|D'_z|=\pi h^2/16$, so $|m_z|\le C\|W/y\|_{L^2(D'_z)}$. Hence

\begin{equation}
|W(z)|\le C\left(\norm{\nabla W}_{L^2(D_z)}+\norm{W/y}_{L^2(D'_z)}\right).
\tag{2.4}
\end{equation}

Lemma 2.1 immediately implies

\[
\|W\|_\infty\le C\|\nabla W\|_2=C\sqrt E,\ \ \ \lim_{|z|\to\infty,\ z\in\overline\HH}W(z)=0.
\]

\qedbox

\result{Lemma 2.3 }
We have

\begin{equation}
D^2W,D^3W,\nabla P,D^2P\in L^2(\HH).
\tag{2.5}
\end{equation}

In particular, $\omega\in H^2(\HH)$ and

\begin{equation}
f(x):=\omega(x,0)\in H^{3/2}(\RR)\subset L^2(\RR).
\tag{2.6}
\end{equation}

\step{Proof. Step 1: second derivatives of velocity.}
Set $F=-(B+W)\cdot\nabla W$ and $M=|B|+\|W\|_\infty$. Then $\|F\|_2\le M\sqrt E$. Write $D_R=B_R(0)\cap\HH$. The scaled local Stokes estimates \cite{Galdi2011}, together with the Poincar\'e inequality at the flat boundary, give

\begin{equation}
\norm{D^2W}_{L^2(D_R)}+\norm{\nabla P}_{L^2(D_R)}
\le C\left(\norm F_{L^2(D_{2R})}+R^{-1}\norm{\nabla W}_{L^2(D_{2R})}\right),
\tag{2.7}
\end{equation}

Here $C$ is independent of $R$. Letting $R\to\infty$ gives $D^2W,\nabla P\in L^2$.

\step{Step 2: one further derivative and the boundary trace.}
Since $\nabla W\in H^1\subset L^4$, we have

\[
\norm{\nabla F}_2\le C\norm{\nabla W}_4^2+M\norm{D^2W}_2<\infty.
\]

The corresponding higher-order Stokes estimate yields

\[
\norm{D^3W}_{L^2(D_R)}+\norm{D^2P}_{L^2(D_R)}
\le C\left(\norm{\nabla F}_{L^2(D_{2R})}+R^{-1}\norm F_{L^2(D_{2R})}
+R^{-2}\norm{\nabla W}_{L^2(D_{2R})}\right).
\]

Letting $R\to\infty$ once more and applying the Sobolev trace theorem, we obtain (2.5) and (2.6). \qedbox

\result{Lemma 2.4 }
After an appropriate choice of the additive constant in $P$,

\begin{equation}
P=\pi+p_s,\qquad p_s=\Py(H_xf),\qquad
\nabla\pi,D^2\pi\in L^2(\HH),\qquad \pi_y(x,0)=0.
\tag{2.8}
\end{equation}

Here $\Py$ is the Poisson extension to the upper half-plane, and $\pi$ is the restriction of the whole-plane quadratic pressure defined below.

\step{Proof. Step 1: reflection and the quadratic tensor.}
Extend $u$ oddly and $v$ evenly across $y=0$, and write $W^e=(u^e,v^e)$. Thus $u^e(x,y)=-u(x,-y)$ and $v^e(x,y)=v(x,-y)$ for $y<0$. The zero trace and incompressibility give

\[
u_x(x,0)=v_x(x,0)=0,\qquad v_y(x,0)=-u_x(x,0)=0.
\]

Then,

\[
\diver W^e=0,\qquad W^e\in L^\infty(\RR^2),\qquad
\nabla W^e,D^2W^e\in L^2(\RR^2).
\]

For $T=W^e\otimes W^e$, the product rule gives

\[
\norm{\nabla T}_2\le C\norm{W^e}_\infty\norm{\nabla W^e}_2,
\quad
\norm{D^2T}_2\le C\left(\norm{W^e}_\infty\norm{D^2W^e}_2+\norm{\nabla W^e}_4^2\right).
\]

Both right-hand sides are finite because $\nabla W^e\in H^1\subset L^4$.

\step{Step 2: the quadratic pressure.}
Define on $\RR^2$

\[
\pi:=\sum_{i,j=1}^2R_iR_jT_{ij},
\]

with the pressure understood modulo an additive constant. The $L^2$ boundedness of the Riesz transforms gives

\[
\norm{\nabla\pi}_2\le C\norm{\nabla T}_2,
\qquad\norm{D^2\pi}_2\le C\norm{D^2T}_2,
\qquad\Delta\pi=-\sum_{i,j=1}^2\partial_i\partial_jT_{ij}.
\]

The diagonal entries of $T$ are even in $y$, whereas its off-diagonal entries are odd. Thus $\pi$ is even. In particular, $\pi_y$ is odd and belongs to $H^1(\RR^2)$, so its boundary trace vanishes.

\step{Step 3: boundary identities for the original pressure.}
Taking the divergence of the momentum equation shows that $\Delta P=-\sum\partial_i\partial_j(W_iW_j)$ in $\HH$. Hence $h=P-\pi$ is harmonic in $\HH$ and $\nabla h\in L^2$. At the boundary, $W=0$, $v_y=0$, and $u_y=-f$. Combining these relations with

\[
\Delta W=(-\omega_y,\omega_x),\qquad
((B+W)\cdot\nabla)W\big|_{y=0}=B_2\partial_yW\big|_{y=0},
\]

we obtain

\begin{equation}
P_x(x,0)=-\omega_y(x,0)+B_2f(x),\qquad P_y(x,0)=f_x(x).
\tag{2.9}
\end{equation}

Since $\pi_y(x,0)=0$, it follows that $h_y(x,0)=f_x$.

\step{Step 4: identification of the harmonic pressure.}
Let $h_0=\Py(H_xf)$. Then

\[
\widehat h_0(k,y)=e^{-|k|y}(-i\operatorname{sgn}k)\widehat f(k),
\qquad \partial_yh_0(x,0)=f_x(x),
\]

and $\|\nabla h_0\|_{L^2(\HH)}^2=\int_\RR|k||\widehat f(k)|^2\dd k<\infty$. The difference $d=h-h_0$ is harmonic, has zero normal derivative on the boundary, and has square-integrable gradient. The even extension of $d$ is harmonic on $\RR^2$ and has square-integrable gradient. Each of its first derivatives is an $L^2$ harmonic function and therefore vanishes. Thus $d$ is constant. Absorbing this constant into the normalization of $P$ proves (2.8).\qedbox

\section[Outflow]{Outflow: $B_2<0$}
Write $B_2=-\gamma$, where $\gamma>0$; no restriction is placed on $B_1$.
We first obtain a normal-direction weight for the vorticity, and then recover $W\in L^2(\HH)$.

\result{Lemma 3.1 }
Set

\[
g(y)=\int_{\RR}\omega^2\dd x,\qquad
d(y)=\int_{\RR}|\nabla\omega|^2\dd x,\qquad
n(y)=\int_{\RR}v\omega^2\dd x.
\]

Then $g\in W^{1,1}(0,\infty)$, and for almost every $y>0$,

\begin{equation}
g'(y)=-\gamma g(y)+n(y)-2\int_y^\infty d(s)\dd s.
\tag{3.2}
\end{equation}

\step{Proof.}
Since $\omega,\omega_y\in L^2$, we have $g'=2\int\omega\omega_y$ and $\int|g'|\le2\|\omega\|_2\|\omega_y\|_2$. Also $g,d,n\in L^1(0,\infty)$, using boundedness of $v$ for the last assertion. Test the vorticity equation $-\Delta\omega+(B+W)\cdot\nabla\omega=0$ with $\eta(y)\omega$, where $\eta\in C_c^\infty(0,\infty)$. Integration by parts gives
\[
\int_0^\infty\eta d+\frac12\int_0^\infty\eta'g'
-\frac12\int_0^\infty\eta'(-\gamma g+n)=0.
\]

Therefore $J=-g'/2-\gamma g/2+n/2$ satisfies $J'=-d$. Since $J,J'\in L^1$, its limit at infinity exists and is zero. Consequently $J(y)=\int_y^\infty d(s)\dd s$, which is (3.2). \qedbox

\result{Lemma 3.2 }
For every $0<\sigma<\gamma$,

\begin{equation}
\int_{\HH}e^{\sigma y}\omega^2<\infty,
\qquad\text{in particular, }y\omega\in L^2(\HH).
\tag{3.3}
\end{equation}

\step{Proof.}
Choose $\varepsilon>0$ so that $\sigma<\gamma-\varepsilon$. Uniform decay gives a height $Y$ such that $\sup_x|v(x,y)|\le\varepsilon$ for $y\ge Y$. Thus $|n|\le\varepsilon g$ there. Dropping the last, nonpositive term in (3.2), we obtain $g'\le-(\gamma-\varepsilon)g$, and hence $g(y)\le g(Y)e^{-(\gamma-\varepsilon)(y-Y)}$. Integration above $y$, together with $\omega\in L^2$ below $Y$, proves the claim.\qedbox

\result{Lemma 3.3.}
Let $W=(u,v)\in H^1_{\rm loc}(\overline\HH;\RR^2)$ satisfy
$\diver W=0$ in $\HH$ and $W|_{y=0}=0$.
Suppose that $\nabla W\in L^2(\HH)$ and $y\omega\in L^2(\HH)$,
where $\omega=v_x-u_y$. Then

\begin{equation}
W\in L^2(\HH),\qquad
\|W\|_2\le2\|y\omega\|_2.
\tag{3.4}
\end{equation}

Moreover,

\[
\|v\|_2\le\|y\omega\|_2,\qquad
\|y\nabla W\|_2^2
=\|y\omega\|_2^2+2\|v\|_2^2
\le3\|y\omega\|_2^2.
\]

\step{Proof. Step 1: horizontal frequency truncation.}
Lemma 2.1 gives $W/y\in L^2(\HH)$ and $\omega\in L^2(\HH)$.
For every $L>0$,

\[
\|W\|_{2,\RR\times(0,L)}
\le L\|W/y\|_2.
\]

Thus $W\in H^1(\RR\times(0,L))$ for every finite $L$, and
$W(\cdot,y)\in L^2(\RR)$ for almost every $y>0$.
We may therefore take its horizontal Fourier transform without assuming
$W\in L^2(\HH)$.

For $\varepsilon>0$, define

\[
\widehat{W_\varepsilon}(k,y)
=\mathbf1_{\{|k|>\varepsilon\}}\widehat W(k,y).
\]

Plancherel's theorem and Tonelli's theorem give

\[
\begin{aligned}
\|W_\varepsilon\|_2^2
&=\int_0^\infty\int_{|k|>\varepsilon}
  |\widehat W(k,y)|^2\dd k\dd y\\
&\le\varepsilon^{-2}
  \int_0^\infty\int_\RR
  |k|^2|\widehat W(k,y)|^2\dd k\dd y
=\varepsilon^{-2}\|\partial_xW\|_2^2<\infty.
\end{aligned}
\]

The frequency projection commutes with derivatives and, on each finite
strip, with the boundary trace. Hence

\[
W_\varepsilon=(u_\varepsilon,v_\varepsilon)\in H_0^1(\HH;\RR^2),
\qquad \diver W_\varepsilon=0.
\]

Its vorticity $\omega_\varepsilon=(v_\varepsilon)_x-(u_\varepsilon)_y$
satisfies

\[
\widehat{\omega_\varepsilon}(k,y)
=\mathbf1_{\{|k|>\varepsilon\}}\widehat\omega(k,y),
\qquad
\|y\omega_\varepsilon\|_2\le \|y\omega\|_2.
\]

The last inequality follows from Plancherel's theorem, since multiplication
by $y$ commutes with the horizontal projection.

\step{Step 2: the stream-function estimate.}
Define $\psi_\varepsilon$ by

\[
\widehat{\psi_\varepsilon}(k,y)
=\begin{cases}
-\widehat{v_\varepsilon}(k,y)/(ik),& |k|>\varepsilon,\\
0,& |k|\le\varepsilon.
\end{cases}
\]

Then $\|\psi_\varepsilon\|_2\le\varepsilon^{-1}\|v_\varepsilon\|_2$.
By incompressibility,

\[
(\psi_\varepsilon)_x=-v_\varepsilon,
\qquad
(\psi_\varepsilon)_y=u_\varepsilon,
\qquad
-\Delta\psi_\varepsilon=\omega_\varepsilon.
\]

The multiplier defining $\psi_\varepsilon$ is bounded on $L_x^2$ and
commutes with the boundary trace. Since $v_\varepsilon|_{y=0}=0$,
we have $\psi_\varepsilon\in H_0^1(\HH)$.

As $\omega_\varepsilon\in L^2(\HH)$, testing the last equation with
$\psi_\varepsilon$ is justified by density and gives

\[
\|W_\varepsilon\|_2^2
=\int_\HH|\nabla\psi_\varepsilon|^2
=\int_\HH\omega_\varepsilon\psi_\varepsilon.
\]

Applying the Hardy inequality used in (2.1a) to $\psi_\varepsilon$, we obtain

\[
\|\psi_\varepsilon/y\|_2
\le2\|(\psi_\varepsilon)_y\|_2
=2\|u_\varepsilon\|_2.
\]

Consequently,

\[
\begin{aligned}
\|u_\varepsilon\|_2^2+\|v_\varepsilon\|_2^2
&\le\|y\omega_\varepsilon\|_2\|\psi_\varepsilon/y\|_2\\
&\le2\|y\omega\|_2\|u_\varepsilon\|_2.
\end{aligned}
\]

Therefore,

\[
\|v_\varepsilon\|_2\le \|y\omega\|_2,
\qquad \|W_\varepsilon\|_2\le2\|y\omega\|_2 .
\]

\step{Step 3: removal of the frequency truncation.}
By monotone convergence and the slice-wise Plancherel identity,

\[
\begin{aligned}
\int_\HH|W|^2
&=\int_0^\infty\int_\RR
  |\widehat W(k,y)|^2\dd k\dd y\\
&=\lim_{\varepsilon\downarrow0}
  \int_0^\infty\int_{|k|>\varepsilon}
  |\widehat W(k,y)|^2\dd k\dd y\\
&\le4\|y\omega\|_2^2.
\end{aligned}
\]

The same argument for $v_\varepsilon$ gives $\|v\|_2\le \|y\omega\|_2$.
This proves (3.4), and in particular $W\in H_0^1(\HH;\RR^2)$.

\step{Step 4: the weighted gradient identity.}
For $R>0$, set

\[
q_R(y)=\frac{R^2y^2}{(R+y)^2}.
\]

For $Z=(z_1,z_2)\in C_c^\infty(\HH;\RR^2)$, integration by parts yields

\[
\begin{aligned}
&\int_\HH q_R\bigl(|\nabla Z|^2-|\curl Z|^2-|\diver Z|^2\bigr)\\
&\quad=2\int_\HH q_R
  \bigl((z_1)_y(z_2)_x-(z_1)_x(z_2)_y\bigr)\\
&\quad=2\int_\HH q_R'z_2(z_1)_x\\
&\quad=2\int_\HH q_R'z_2\diver Z
  +\int_\HH q_R''z_2^2.
\end{aligned}
\]

For fixed $R$, the functions $q_R$, $q_R'$, and $q_R''$ are bounded.
We may therefore pass to $W$ by approximation in $H_0^1(\HH;\RR^2)$.
Since $\diver W=0$, this gives

\[
\int_\HH q_R|\nabla W|^2
=\int_\HH q_R\omega^2+\int_\HH q_R''v^2.
\]

Finally,

\[
q_R\uparrow y^2,
\qquad
q_R''(y)=\frac{2R^3(R-2y)}{(R+y)^4}\longrightarrow2,
\qquad |q_R''|\le2.
\]

Since $y\omega,v\in L^2(\HH)$, monotone convergence theorem for the first two
integrals and dominated convergence theorem for the last one yield

\[
\|y\nabla W\|_2^2
=\|y\omega\|_2^2+2\|v\|_2^2
\le3\|y\omega\|_2^2.
\]
\qedbox

\result{Theorem 3.4 (Rigidity in the outflow case).}
Under the assumptions of Theorem 1.1, if $B_2<0$, then $U\equiv B$ in $\HH$ and $P$ is constant.

\step{Proof. Step 1: the energy identity.}
Lemmas 3.2 and 3.3 imply $W\in L^2$, hence $W\in H_0^1(\HH)$. We also have $W\in L^\infty$ and $\nabla P\in L^2$. Choose a standard radial cutoff $\eta_R$ and write

\[
S_R=(B_{2R}\setminus B_R)\cap\HH,\qquad P_R=|S_R|^{-1}\int_{S_R}P.
\]

Testing the momentum equation with $W\eta_R^2$, incompressibility and the zero trace give

\[
\int((B+W)\cdot\nabla W)\cdot W\eta_R^2
=-\frac12\int|W|^2(B+W)\cdot\nabla(\eta_R^2).
\]

 Rearranging yields

\[
\begin{aligned}
\int\eta_R^2|\nabla W|^2
={}&-2\int\eta_R\sum_{j=1}^2(\partial_jW\cdot W)\partial_j\eta_R\\
&+\frac12\int|W|^2(B+W)\cdot\nabla(\eta_R^2)
+\int(P-P_R)W\cdot\nabla(\eta_R^2).
\end{aligned}
\]

\step{Step 2: passage to infinity.}
Denote the three terms on the right by $I_R,J_R,K_R$. The Poincar\'e inequality on the half-annulus gives $\|P-P_R\|_{2,S_R}\le CR\|\nabla P\|_{2,S_R}$, and therefore

\[
\begin{aligned}
|I_R|&\le CR^{-1}\norm{\nabla W}_{2,S_R}\norm W_{2,S_R},\\
|J_R|&\le CR^{-1}(|B|+\norm W_\infty)\norm W_{2,S_R}^2,\\
|K_R|&\le C\norm{\nabla P}_{2,S_R}\norm W_{2,S_R}.
\end{aligned}
\]

All three terms tend to zero because $W,\nabla W,\nabla P\in L^2$. Letting $R\to\infty$ gives $E=0$. The zero boundary trace then implies $W=0$, and hence $U=B$ and $P$ is constant. \qedbox

% section4
\section[Inflow]{Inflow: $B_2>0$}
The tangential component $B_1$ is arbitrary. We recall the results from Section 2:

\[
P=\pi+p_s,\qquad
\pi=\sum_{i,k=1}^2R_iR_k(W_i^eW_k^e),\qquad
p_s=\PP_y(H_xf),\quad f=\omega(\cdot,0),
\]

where $u^e$ is odd and $v^e$ is even in $y$.  Introduce the whole-plane notation

\begin{equation}
\begin{gathered}
t=\rho+|y|,\qquad \chi_\rho=\frac\rho t,\qquad
\vartheta_\rho=\chi_\rho^2,\qquad w_\rho=\frac{\vartheta_\rho}\rho,\qquad
A_\rho=\int_\HH w_\rho|W|^2.
\end{gathered}\tag{4.1}
\end{equation}

 The Hardy inequality gives $A_\rho\le\rho\|W/y\|_2^2<\infty$ and $\|\chi_\rho W\|_2^2=\rho A_\rho$. Since $\nabla(\chi_\rho W)=\chi_\rho\nabla W+\chi_\rho'W\otimes e_2$, it follows that $\chi_\rho W\in H_0^1(\HH)$. For $y>0$,

\begin{equation}
\chi_\rho'=-w_\rho,\qquad \chi_\rho''=\frac{2\rho}{t^3},\qquad
w_\rho'=-\frac{2w_\rho}{t}.\tag{4.2}
\end{equation}

\subsection{The quadratic pressure term}
\result{Lemma 4.1 }
Let $\pi$ be the pressure generated by $W^e\otimes W^e$ in Lemma 2.4. Then, for every $\rho>0$,

\begin{equation}
J_{\pi,\rho}:=-\int_\HH\nabla\pi\cdot(\chi_\rho W),\qquad
|J_{\pi,\rho}|\le C\|W\|_\infty A_\rho.\tag{4.3}
\end{equation}

The constant is independent of $\rho$ and $W$.

\step{Proof. Step 1: a pressure identity.}
From the definition of $\pi$, we obtain

\[
\pi=R_x^2\bigl((u^e)^2\bigr)+2R_xR_y(u^ev^e)+R_y^2\bigl((v^e)^2\bigr).
\]

 Define

\[
T=\TT(W):=R_xR_y\bigl((u^e)^2-(v^e)^2\bigr)-2R_x^2(u^ev^e).
\]

 Using $R_x^2+R_y^2=-I$, $\partial_xR_y^2=\partial_yR_xR_y$, and $\partial_xR_xR_y=\partial_yR_x^2$, we obtain

\begin{equation}
\begin{aligned}
\partial_x\pi
&=-\partial_x\bigl((u^e)^2\bigr)-\partial_yR_xR_y\bigl((u^e)^2\bigr)\\
&\quad+2\partial_yR_x^2(u^ev^e)+\partial_yR_xR_y\bigl((v^e)^2\bigr)\\
&=-\partial_x\bigl((u^e)^2\bigr)-\partial_yT\qquad\text{in }\RR^2.
\end{aligned}\tag{4.4}
\end{equation}

\step{Step 2: the weighted estimate for the odd output.}
Recall that $w_\rho=\chi_\rho^2/\rho$. The functions $(u^e)^2$ and $(v^e)^2$ are even, and $u^ev^e$ is odd. Applying Lemma B.3 to the $R_xR_y$ terms and the $R_x^2$ term therefore gives

\begin{equation}
\begin{aligned}
\|\sqrt{w_\rho}T\|_{2,\HH}
&\le C\|\sqrt{w_\rho}(W^e\otimes W^e)\|_{2,\RR^2}\\
&\le C\|W\|_\infty A_\rho^{1/2}.
\end{aligned}\tag{4.5}
\end{equation}

\step{Step 3: a weighted stream-function estimate.}
 Using $\partial_y(t^{-3})=\partial_y{(\rho+|y|)^{-3}}=-3t^{-4}, y>0$ and $\psi(x,0)=0$, integration by parts yields

\[
\begin{aligned}
\int_{\RR}\int_0^L\frac{\rho\psi^2}{t^4}\,dy\,dx
&=-\frac{\rho}{3(\rho+L)^3}\int_{\RR}|\psi(x,L)|^2\,dx
 +\frac23\int_{\RR}\int_0^L\frac{\rho\psi u^e}{t^3}\,dy\,dx\\
&\le\frac23\int_{\RR}\int_0^L\frac{\rho|\psi u^e|}{t^3}\,dy\,dx.
\end{aligned}
\]

 By Cauchy--Schwarz, we obtain

\[
\begin{aligned}
\int_{\RR}\int_0^L\frac{\rho\psi^2}{t^4}\,dy\,dx
&\le\frac23
\left(\int_{\RR}\int_0^L\frac{\rho\psi^2}{t^4}\,dy\,dx\right)^{1/2}
\left(\int_{\RR}\int_0^L\frac{\rho|u^e|^2}{t^2}\,dy\,dx\right)^{1/2}.
\end{aligned}
\]

 Since $w_\rho=\rho/t^2$, we obtain

\[
\int_{\RR}\int_0^L\frac{\rho\psi^2}{t^4}\,dy\,dx
\le\frac49\int_{\RR}\int_0^L\frac{\rho|u^e|^2}{t^2}\,dy\,dx
\le\frac49\|\sqrt{w_\rho}u^e\|_{2,\HH}^2.
\]

 Letting $L\to\infty$ and using monotone convergence gives

\begin{equation}
\begin{aligned}
\left\|\sqrt{w_\rho}\frac\psi t\right\|_2
&=\left(\int_{\RR}\int_0^\infty\frac{\rho\psi^2}{t^4}\,dy\,dx\right)^{1/2}\\
&\le\frac23\|\sqrt{w_\rho}u^e\|_{2,\HH}\le\frac23 A_\rho^{1/2}.
\end{aligned}\tag{4.6}
\end{equation}

\step{Step 4: the pressure-work identity.}
Choose a smooth even cutoff $\sigma_N(y)$ equal to one for $|y|\le N$ and zero for $|y|\ge2N$, with $|\sigma_N'|\le C/N$. Set $F_N=\sigma_NW^e\otimes W^e$ and

\[
\pi_N=\sum_{i,k=1}^2R_iR_k(F_N)_{ik},\qquad
T_N=R_xR_y\bigl((F_N)_{11}-(F_N)_{22}\bigr)-2R_x^2(F_N)_{12}.
\]

Then $\pi_N,T_N\in H^1(\RR^2)$, and the analogue of (4.4) is $\partial_x\pi_N=-\partial_x(F_N)_{11}-\partial_yT_N$. Since $\diver(\chi_\rho W)=\chi_\rho'v^e=-w_\rho v^e$ and $v^e=-\psi_x$, successive integrations by parts give

\[
\begin{aligned}
-\int_\HH\nabla\pi_N\cdot(\chi_\rho W)
&=\int_\HH\pi_N\operatorname{div}(\chi_\rho W)
=-\int_\HH w_\rho\pi_Nv^e\\
&=\int_\HH w_\rho\pi_N\psi_x
=-\int_\HH w_\rho\psi\,\partial_x\pi_N\\
&=\int_\HH w_\rho\psi\,\partial_x(F_N)_{11}
 +\int_\HH w_\rho\psi\,\partial_yT_N\\
&=-\int_\HH w_\rho\psi_x(F_N)_{11}
 -\int_\HH T_N\partial_y(w_\rho\psi)\\
&=\int_\HH w_\rho(F_N)_{11}v^e
 -\int_\HH T_N(w_\rho u^e+w_\rho'\psi).
\end{aligned}
\]

Here $w_\rho$ is independent of $x$, and $\partial_y(w_\rho\psi)=w_\rho u^e+w_\rho'\psi,\chi_\rho W,w_\rho\psi\in H_0^1(\HH)$ . Passing to the limit as explained in Step 5, we obtain

\begin{equation}
J_{\pi,\rho}
=\int_\HH w_\rho (u^e)^2v^e
-\int_\HH T\bigl(w_\rho u^e+w_\rho'\psi\bigr).\tag{4.7}
\end{equation}

 Consequently,

\[
\begin{aligned}
|J_{\pi,\rho}|
&\le\|W\|_\infty A_\rho
+\|\sqrt{w_\rho}T\|_2
\left(\|\sqrt{w_\rho}u^e\|_{2,\HH}+2\left\|\sqrt{w_\rho}\frac\psi t\right\|_2\right)\\
&\le C\|W\|_\infty A_\rho.
\end{aligned}
\]

It remains to justify the limiting procedure.

\step{Step 5: removal of the cutoff.}
Write $F=W^e\otimes W^e$. Section 2 and the Hardy inequality imply $\nabla F,F/y\in L^2$; also $F$ is square integrable on every strip of finite height. Hence $F_N\in H^1$, and

\[
\|\nabla(F_N-F)\|_2
\le\|1_{\{|y|>N\}}\nabla F\|_2
+C\|1_{\{|y|>N\}}F/y\|_2\longrightarrow0.
\]

It follows that $\nabla\pi_N\to\nabla\pi$ in $L^2$. Also $\chi_\rho F_N\to \chi_\rho F$ in $L^2$, so Lemma B.3 gives $\chi_\rho T_N\to \chi_\rho T$ in $L^2$.

By (4.6) and $\psi_x=-v^e$, $\psi_y=u^e$, we have $w_\rho\psi\in H_0^1(\HH)$. Thus all expressions in Step 4 are legitimate Sobolev pairings. The preceding convergences, together with $\chi_\rho W\in L^2$, allow passage to the limit in each term and prove (4.7).\qedbox

\result{Remark }
The next lemma uses the same calculation for a truncated source. Let $0\le\theta\le1$ be smooth and even in $y$, with compactly supported gradient, and set

\[
F_\theta=\theta W^e\otimes W^e,\qquad
\pi_\theta=\sum_{i,k}R_iR_k(F_\theta)_{ik},
\]

\[
T_\theta=R_xR_y\bigl((F_\theta)_{11}-(F_\theta)_{22}\bigr)
-R_x^2\bigl((F_\theta)_{12}+(F_\theta)_{21}\bigr).
\]

Then $\partial_x\pi_\theta=-\partial_x(F_\theta)_{11}-\partial_yT_\theta$, and $T_\theta$ is odd. Steps 2--5 give

\[
\begin{aligned}
-\int_\HH\nabla\pi_\theta\cdot(\chi_\rho W)
&=\int_\HH w_\rho(F_\theta)_{11}v^e
-\int_\HH T_\theta(w_\rho u^e+w_\rho'\psi),\\
\left|-\int_\HH\nabla\pi_\theta\cdot(\chi_\rho W)\right|
&\le C\|\sqrt{w_\rho}F_\theta\|_{2,\RR^2}\,A_\rho^{1/2}.
\end{aligned}
\]

\result{Lemma 4.2 }
For every $\varepsilon>0$,

\begin{equation}
|J_{\pi,\rho}|\le C\varepsilon A_\rho+o_\rho(1),\qquad\rho\to\infty.\tag{4.8}
\end{equation}

where $C$ is independent of $\varepsilon$ and $\rho$. The remainder tends to zero with $\varepsilon$ fixed.

\step{Proof. Step 1: a compact part and a tail.}
Since $W^e\in C_0(\RR^2)$, choose $R_\varepsilon$ so that $|W^e|\le\varepsilon$ for $|z|\ge R_\varepsilon$. Let $0\le\eta_\varepsilon\le1$ be a radial cutoff, equal to one on $B_{R_\varepsilon}$ and zero outside $B_{2R_\varepsilon}$, and set

\[
F_c=\eta_\varepsilon W^e\otimes W^e,\qquad
F_t=(1-\eta_\varepsilon)W^e\otimes W^e.
\]

Then $\nabla\pi=\nabla\pi_c+\nabla\pi_t$, with $\pi_c\in H^1(\RR^2)$. Since $|F_t|\le\varepsilon|W^e|$, the preceding remark gives

\[
\left|-\int_\HH\nabla\pi_t\cdot(\chi_\rho W)\right|
\le C\|\sqrt{w_\rho}F_t\|_2 A_\rho^{1/2}
\le C\varepsilon A_\rho.
\]

\step{Step 2: the compact contribution.}
Using $\chi_\rho W\in H_0^1$ and incompressibility,

\[
-\int_\HH\nabla\pi_c\cdot(\chi_\rho W)=\int_\HH\pi_c\chi_\rho'v.
\]

Moreover, $|y\chi_\rho'|=\rho y/(\rho+y)^2\le1/4$, and this quantity tends to zero for each fixed $y$ as $\rho\to\infty$. Since $v/y\in L^2$, dominated convergence yields

\[
\|\chi_\rho'v\|_2^2=\int_\HH|y\chi_\rho'|^2|v/y|^2\longrightarrow0.
\]

For fixed $\varepsilon$, the function $\pi_c$ is independent of $\rho$. Its pressure work is therefore bounded by $\|\pi_c\|_2\|\chi_\rho'v\|_2\to0$. Combining the two contributions proves (4.8).\qedbox

\subsection{The weighted energy identity and rigidity}
\result{Theorem 4.3 (Rigidity in the inflow case).}
Under the assumptions of Theorem 1.1, if $B_2>0$, then $U\equiv B$ in $\HH$ and $P$ is constant.

\step{Proof. Step 1: the weighted energy identity.}
Fix $\rho>0$. Since $\chi_\rho W\in H_0^1$, $W$ is bounded, and $\nabla W,\nabla P\in L^2$, testing the momentum equation with $\chi_\rho W$ gives

\[
\int_\HH\nabla W:\nabla(\chi_\rho W)
+\int_\HH(B\cdot\nabla W)\cdot(\chi_\rho W)
+\int_\HH(W\cdot\nabla W)\cdot(\chi_\rho W)
+\int_\HH\nabla P\cdot(\chi_\rho W)=0.
\]

For brevity, write

\[
\begin{aligned}
D_\rho&=\tfrac12\int_\HH\chi_\rho''|W|^2,&
N_\rho&=\tfrac12\int_\HH\chi_\rho'v|W|^2,\\
J_{\pi,\rho}&=-\int_\HH\nabla\pi\cdot(\chi_\rho W),&
J_{s,\rho}&=\int_\HH\chi_\rho'p_sv.
\end{aligned}
\]

Integration by parts in each term gives

\[
\begin{aligned}
\int_\HH\nabla W:\nabla(\chi_\rho W)
&=\int_\HH\chi_\rho|\nabla W|^2+\tfrac12\int_\HH\chi_\rho'\partial_y|W|^2\\
&=\int_\HH\chi_\rho|\nabla W|^2-D_\rho,\\
\int_\HH(B\cdot\nabla W)\cdot(\chi_\rho W)
&=-\tfrac{B_2}{2}\int_\HH\chi_\rho'|W|^2=\tfrac{B_2}{2}A_\rho,\\
\int_\HH(W\cdot\nabla W)\cdot(\chi_\rho W)
&=-\tfrac12\int_\HH\chi_\rho'v|W|^2=-N_\rho,\\
\int_\HH\nabla P\cdot(\chi_\rho W)&=-J_{\pi,\rho}-J_{s,\rho}.
\end{aligned}
\]

Substitution into the test identity yields

\begin{equation}
\int_\HH\chi_\rho|\nabla W|^2+\frac{B_2}{2}A_\rho
=D_\rho+N_\rho+J_{\pi,\rho}+J_{s,\rho}.\tag{4.9}
\end{equation}

\step{Step 2: the viscous error.}
By (4.2) and the dominated convergence theorem, we obtain

\[
D_\rho=\int_\HH\frac{\rho y^2}{(\rho+y)^3}\left|\frac Wy\right|^2\longrightarrow0.
\]

\step{Step 3: the nonlinear terms.}
Fix $\varepsilon>0$ and choose $R_\varepsilon$ so that $|W|\le\varepsilon$ outside $B_{R_\varepsilon}$. Since $|\chi_\rho'|=w_\rho\le1/\rho$,

\[
|N_\rho|\le\frac1{2\rho}\int_{B_{R_\varepsilon}\cap\HH}|v||W|^2
+\frac\varepsilon2\int_{\HH\setminus B_{R_\varepsilon}}w_\rho|W|^2.
\]

Together with Lemma 4.2, this gives, for fixed $\varepsilon$,

\begin{equation}
|N_\rho|\le\tfrac\varepsilon2 A_\rho+o_\rho(1),\qquad
|J_{\pi,\rho}|\le C\varepsilon A_\rho+o_\rho(1).\tag{4.10}
\end{equation}

\step{Step 4: the harmonic pressure.}
Let $g=H_xf\in L^2(\RR)$. Since $p_s=\Py g$ and the Poisson semigroup is contractive on $L^2$,

\[
S_\rho:=\int_\HH w_\rho|p_s|^2
\le\|g\|_2^2\int_0^\infty\frac\rho{(\rho+y)^2}\,dy=\|g\|_2^2<\infty.
\]

The substitution $y=\rho s$ gives

\[
S_\rho=\int_0^\infty(1+s)^{-2}\|\PP_{\rho s}g\|_{L^2(\RR)}^2\,ds.
\]

For every $s>0$, Plancherel's theorem and dominated convergence theorem imply

\[
\|\PP_{\rho s}g\|_2^2=\int_\RR e^{-2\rho s|k|}|\widehat g(k)|^2\,dk
\longrightarrow0\qquad(\rho\to\infty).
\]

Also $\|\PP_{\rho s}g\|_2\le\|g\|_2$, and $(1+s)^{-2}$ is integrable on $(0,\infty)$. Hence $S_\rho\to0$. Cauchy--Schwarz and Young's inequality give

\begin{equation}
|J_{s,\rho}|\le\left(\int_\HH w_\rho v^2\right)^{1/2}S_\rho^{1/2}
\le\frac{B_2}{8}A_\rho+\frac2{B_2}S_\rho.\tag{4.11}
\end{equation}

\step{Step 5: absorption and passage to the limit.}
Since $B_2>0$, choose $\varepsilon>0$ sufficiently small that

\[
\frac\varepsilon2+C\varepsilon\le\frac{B_2}{8}.
\]

Insert (4.10)--(4.11) into (4.9) and absorb the $A_\rho$ terms into the left-hand side. We obtain

\[
\int_\HH\chi_\rho|\nabla W|^2+\frac{B_2}{4}A_\rho
\le D_\rho+o_\rho(1)+\frac2{B_2}S_\rho\longrightarrow0.
\]

The quantity $A_\rho$ is finite for each fixed $\rho$ by (4.1), so this absorption is legitimate. Finally, $0\le\chi_\rho\le1$, $\chi_\rho\to1$, and the Dirichlet integral is finite. Therefore

\[
\int_\HH|\nabla W|^2=0.
\]

Thus $W$ is constant, and its zero trace gives $W=0$. Then $U\equiv B$ in $\HH$ and $P$ is constant.\qedbox

% section5
\section{Tangential backgrounds: identities and weighted estimates}
In this section and the next, $B=(b,0)$ with $b>0$; a horizontal reflection reduces $b<0$ to this case. We first prove $v\in L^2(\HH)$ and then use a truncated stream function to establish rigidity. Tangential transport produces no favorable normal term in the usual energy identity. The normal momentum equation, however, contains $bv_x$. We therefore express its remaining terms as horizontal derivatives and use the other momentum equation and the boundary conditions to obtain an identity for $bv$ itself. This reduces the problem to a quadratic bulk term and a boundary-vorticity term.

\subsection{An identity from the momentum equations}
For now, every quadratic source comes from the actual solution. Put $f=\omega(\cdot,0)$ and $h_f=\Py f$. Lemma 2.4 gives

\[
P=\pi+p_s,\qquad p_s=\Py(H_xf),\qquad
\pi=R_x^2(u^2)+2R_xR_y(uv)+R_y^2(v^2).
\]

The products are extended to the whole plane according to the odd extension of $u$ and even extension of $v$. To combine the normal convection term with the quadratic pressure, define

\[
G(W\otimes W):=R_xR_y(u^2-v^2)+(R_y^2-R_x^2)(uv).
\]

where $G(W\otimes W)$ is odd in $y$.

\result{Lemma 5.1 (An identity eliminating the quadratic pressure).}
The solution satisfies

\begin{equation}
bv+G(W\otimes W)=\omega-h_f\qquad\text{in }\HH.
\tag{5.1}\label{eq:tangent-identity}
\end{equation}

\step{Proof. Step 1: the normal momentum equation.}
Using $\diver W=0$ and $\Delta v=\omega_x$, we obtain

\[
bv_x=\omega_x-\partial_x(uv)-\partial_y(v^2)-P_y.
\]

We seek a horizontal derivative representation of $\partial_x(uv)+\partial_y(v^2)+\pi_y$. The Riesz multipliers satisfy

\[
R_x^2+R_y^2=-I,\qquad
\partial_yR_xR_y=\partial_xR_y^2,\qquad
\partial_yR_x^2=\partial_xR_xR_y.
\]

Differentiating the expression for $\pi$ term by term gives

\[
\begin{aligned}
\partial_y{\pi}
&=\partial_y\{R_x^2(u^2)+2R_xR_y(uv)-R_x^2(v^2)-v^2\}
\\
&=\partial_xR_xR_y(u^2-v^2)+\partial_x\{R_y^2-I-R_x^2\}(uv)-\partial_y(v^2)\\
&=\partial_x\{G(W\otimes W)\}-\partial_x(uv)-\partial_y(v^2).
\end{aligned}
\]

 we find

\[
\partial_x(uv)+\partial_y(v^2)+\pi_y=\partial_x\{G(W\otimes W)\}.
\]

Therefore
\[
bv_x=\omega_x-\partial_xG(W\otimes W)-(p_s)_y.
\]

\step{Step 2: the harmonic pressure.}
Since the symbol of $H_x$ is $-i\operatorname{sgn}k$,

\[
\widehat{(p_s)_y}(k,y)
=-|k|e^{-|k|y}(-i\operatorname{sgn}k)\widehat f(k)
=ik e^{-|k|y}\widehat f(k).
\]

Therefore $(p_s)_y=\partial_xh_f$, and

\[
\partial_x\{bv+G(W\otimes W)-\omega+h_f\}=0.
\]

\step{Step 3: determination of the constant.}
The same multiplier identities give, for the horizontal pressure combination,

\[
\begin{aligned}
\pi_x+\partial_x(u^2)+\partial_y(uv)
&=-\partial_xR_y^2(u^2-v^2)
  +\partial_y(R_x^2-R_y^2)(uv)\\
&=-\partial_yG(W\otimes W).
\end{aligned}
\]

Since $\Delta u=-\omega_y$, the horizontal momentum equation becomes

\[
bu_x=-\omega_y+\partial_yG(W\otimes W)-(p_s)_x.
\]

The other derivative of the harmonic pressure satisfies

\[
\widehat{(p_s)_x}(k,y)
=ik(-i\operatorname{sgn}k)e^{-|k|y}\widehat f(k)
=|k|e^{-|k|y}\widehat f(k)
=-\widehat{(h_f)_y}(k,y).
\]

Combining these identities with $u_x+v_y=0$ gives

\[
\partial_y\{bv+G(W\otimes W)-\omega+h_f\}
=b(v_y+u_x)=0.
\]

 By Lemma 2.4, $\nabla G\in L^2$, and the odd representative of $G$ has zero boundary trace. Also $v(\cdot,0)=0$ and $\omega(\cdot,0)=h_f(\cdot,0)=f$. These boundary values prove \eqref{eq:tangent-identity}.\qedbox

The vorticity term in \eqref{eq:tangent-identity} is already in $L^2$. To estimate $v$, it remains to control $G(W\otimes W)$ and $h_f$. The former is a quadratic bulk operator. The latter contains the actual boundary vorticity, and $f\in L^2(\RR)$ alone is insufficient, because

\[
\norm{\Py f}_{L^2(\HH)}^2
=\int_{\RR}\frac{|\widehat f(k)|^2}{2|k|}\dd k
\]

may diverge at $k=0$. We must use the equation and both zero-velocity boundary conditions to determine the relation between $f$ and the quadratic source.

\subsection{The boundary response of the solution}
Keep $F=W\otimes W$ for the moment, and rewrite the  equation as

\[
-\Delta W+b\partial_xW+\nabla P=-\diver F.
\]

\result{Lemma 5.2 (Boundary vorticity generated by the quadratic source).}
Let $F=W\otimes W$ and $f=\omega(\cdot,0)$. Then

\begin{equation}
f=T_b(W\otimes W),\qquad
\widehat{T_b(F)}(k)
=\frac{\mathcal T_\lambda S(F)(k)-\mathcal T_{|k|}S(F)(k)}{\lambda-{|k|}}.
\tag{5.2}\label{eq:actual-wall}
\end{equation}

The quantities $S(F)$, $\lambda$, and $\mathcal T_z$ are introduced in the proof through the vorticity equation and the boundary conditions on the stream function.

\step{Proof. Step 1: the scalar Oseen equation.}
Define the vorticity source

\[
S(F):=-\curl\diver F
=(\partial_y^2-\partial_x^2)F_{12}
+\partial_x\partial_y(F_{11}-F_{22}).
\]

Taking the curl eliminates the pressure and gives

\begin{equation}
-\Delta\omega+b\partial_x\omega=S(F).
\tag{5.3}\label{eq:oseen-vorticity}
\end{equation}

After Fourier transformation in $x$, set, for $k\ne0$,

\[
\lambda=\sqrt{{|k|}^2+ibk},\qquad\Re\lambda>0.
\]

Then

\begin{equation}
(-\partial_y^2+\lambda^2)\widehat\omega=\widehat{S(F)}.
\tag{5.4}\label{eq:scalar-ode}
\end{equation}

Since $\omega\in L^2(\HH)$, we have
$\widehat\omega(k,\cdot)\in L^2(0,\infty)$ for almost every
$k\ne0$. This excludes the growing homogeneous mode
$e^{\lambda y}$. The coefficient of the remaining mode
$e^{-\lambda y}$ is determined by the compatibility condition
derived below from the zero-velocity boundary conditions.

\step{Step 2: the constraint from the  stream function.}
Introduce the Laplace-type integral in the normal variable,

\[
\mathcal T_zH(k):=\int_0^\infty e^{-z(k)y}\widehat H(k,y)\dd y,
\qquad z={|k|},\lambda.
\]

 Since $\widehat\omega=(-\partial_y^2+|k|^2)\widehat\psi$, integration by parts gives

\begin{equation}
\mathcal T_{|k|}\omega
=\partial_y\widehat\psi(k,0)+{|k|}\widehat\psi(k,0)=0,
\tag{5.5}\label{eq:clamped-moment}
\end{equation}

where we use both $\psi=0$ and $\psi_y=0$ at the boundary. It is the vorticity compatibility condition imposed by the full zero-velocity trace.

\step{Step 3: identification of the boundary value.}
Let

\[
D_\lambda(y,z)=\frac{e^{-\lambda|y-z|}-e^{-\lambda(y+z)}}{2\lambda}
\]

be the Dirichlet Green kernel on the half-line. The decaying representation of \eqref{eq:scalar-ode} is

\begin{equation}
\widehat\omega(k,y)=e^{-\lambda y}\widehat f(k)
+\int_0^\infty D_\lambda(y,z)\widehat{S(F)}(k,z)\dd z.
\tag{5.6}\label{eq:green-vorticity}
\end{equation}

 Direct integration of the kernel gives

\[
\int_0^\infty e^{-{|k|}y}D_\lambda(y,z)\dd y
=\frac{e^{-{|k|}z}-e^{-\lambda z}}{\lambda^2-{|k|}^2}.
\]

Insert \eqref{eq:green-vorticity} into \eqref{eq:clamped-moment} to obtain

\[
0=\frac{\widehat f}{{|k|}+\lambda}
+\frac{\mathcal T_{|k|}S(F)-\mathcal T_\lambda S(F)}{\lambda^2-{|k|}^2}.
\]

Using $\lambda^2-|k|^2=(\lambda-|k|)(\lambda+|k|)$ gives \eqref{eq:actual-wall}.
\qedbox

\subsection{Tensor sources and the quadratic expansion}
We now extend the notation to general symmetric tensors. This allows the decomposition $W\otimes W=A\otimes A+K$ to be used linearly at the level of the source, with separate estimates for a small tail and a fixed compact part.

When applying whole-plane Riesz transforms to a half-plane tensor, extend its diagonal entries evenly and its off-diagonal entries oddly, then restrict the output to $\HH$. Define

\begin{equation}
\begin{aligned}
G(F)&=R_xR_y(F_{11}-F_{22})+(R_y^2-R_x^2)F_{12},\\
S(F)&=(\partial_y^2-\partial_x^2)F_{12}
       +\partial_x\partial_y(F_{11}-F_{22}).
\end{aligned}
\tag{5.7}\label{eq:general-GS}
\end{equation}

The operator $G$ combines the normal quadratic term with the pressure and has odd output. The operator $S=-\curl\diver$ produces the vorticity source. For smooth compactly supported sources, define

\begin{equation}
\mathcal T_zH(k)=\int_0^\infty e^{-z(k)y}\widehat H(k,y)\dd y,
\qquad
\widehat{T_b(F)}(k)=\frac{\mathcal T_\lambda S(F)-\mathcal T_{|k|}S(F)}{\lambda-{|k|}},
\quad k\ne0.
\tag{5.8}\label{eq:general-T}
\end{equation}

where $\lambda=\sqrt{|k|^2+ibk}$ and $\Re\lambda>0$. The boundary response $T_b$ is linear in the tensor source. If $Z_F$ is the zero-boundary Oseen response to $F$, the calculation in Lemma 5.2 gives $T_b(F)=\curl Z_F(\cdot,0)$. A general input field $A$ need not solve that response equation. In particular, $T_b(A\otimes A)$ cannot in general be identified with $\curl A(\cdot,0)$; the identity for the actual solution is supplied separately by Lemma 5.2.

Let $\Pil$ denote the orthogonal projection onto $0<|k|\le b/4$. On this interval,

\begin{equation}
|\lambda|\asymp\sqrt{b{|k|}},\qquad
|\lambda-{|k|}|\asymp\sqrt{b{|k|}},\qquad
|\lambda+{|k|}|\le Cb,\qquad
(\lambda-{|k|})(\lambda+{|k|})=ibk.
\tag{5.9}\label{eq:low-frequency}
\end{equation}

Only this low-frequency portion of the boundary response requires a refined estimate. At high frequencies,

\[
\norm{\Py(1-\Pil)f}_2^2
=\int_{|k|>b/4}\frac{|\widehat f(k)|^2}{2|k|}\dd k
\le Cb^{-1}\norm f_2^2.
\]

The operators $\Py$ and $\mathcal T_z$ act in opposite directions: $\Py$ extends boundary data harmonically into the half-plane, whereas $\mathcal T_z$ integrates a bulk source in the normal variable.

We record the quadratic expansion for $A=(r,s)$. First assume that $A$ is smooth, compactly supported, divergence-free, and zero on the boundary. Put $C_0=rs$ and $D=r^2-s^2$. In horizontal Fourier variables,

\[
\widehat{S(A\otimes A)}=(\partial_y^2+{|k|}^2)\widehat C_0+ik\partial_y\widehat D.
\]

Since $D=C_0=\partial_yC_0=0$ on the boundary, integration by parts gives, for $z=|k|,\lambda$,

\[
\begin{aligned}
\mathcal T_z(\partial_yC_0)&=z\mathcal T_zC_0,\\
\mathcal T_z(\partial_y^2C_0)&=z^2\mathcal T_zC_0,\\
\mathcal T_z(\partial_yD)&=z\mathcal T_zD.
\end{aligned}
\]

Hence $\mathcal T_zS(A\otimes A)=(z^2+|k|^2)\mathcal T_zC_0+ikz\mathcal T_zD$. Substituting into \eqref{eq:general-T}, we obtain

\begin{equation}
\widehat{T_b(A\otimes A)}
=\frac{(\lambda^2+{|k|}^2)\mathcal T_\lambda C_0-2{|k|}^2\mathcal T_{|k|}C_0
+ik\lambda\mathcal T_\lambda D-ik{|k|}\mathcal T_{|k|}D}{\lambda-{|k|}}.
\tag{5.10a}\label{eq:wall-expansion}
\end{equation}

Incompressibility, $r_x=-s_y$, allows one horizontal derivative of the tangential square to be transferred to the normal component. Since $rs$ vanishes at the boundary, for $z=|k|,\lambda$,

\begin{equation}
ik\mathcal T_z(r^2)=-2\mathcal T_z(rs_y)
=2\mathcal T_z(r_ys)-2z\mathcal T_z(rs).
\tag{5.10b}\label{eq:transfer-square}
\end{equation}

Apply this identity first to the harmonic $r^2$ term in \eqref{eq:wall-expansion}. The resulting $2|k|^2\Tk(rs)$ cancels the original harmonic mixed term, leaving

\begin{equation}
\begin{aligned}
\widehat{T_b(A\otimes A)}
={}&\frac{ik\lambda\Tl(r^2)+(\lambda^2+|k|^2)\Tl(rs)-2|k|\Tk(r_ys)}{\lambda-|k|}\\
&+ik\frac{|k|\Tk(s^2)-\lambda\Tl(s^2)}{\lambda-|k|}.
\end{aligned}
\tag{5.10}\label{eq:wall-reduced}
\end{equation}

This is the form used for the finite-input estimate below. The harmonic tangential square has been replaced by a product containing $s$, with its full low-frequency coefficient retained. Applying \eqref{eq:transfer-square} also to the Oseen square, and using $(\lambda^2+|k|^2)-2\lambda^2=-(\lambda-|k|)(\lambda+|k|)$, gives

\begin{equation}
\begin{aligned}
\widehat{T_b(A\otimes A)}
={}&\frac{2\lambda\Tl(r_ys)-2|k|\Tk(r_ys)}{\lambda-|k|}
-(\lambda+|k|)\Tl(rs)\\
&+ik\frac{|k|\Tk(s^2)-\lambda\Tl(s^2)}{\lambda-|k|}.
\end{aligned}
\tag{5.10c}\label{eq:wall-cancel}
\end{equation}

All three sources in this formula, $r_ys$, $rs$, and $s^2$, contain the normal component. This form is suited to estimates on strips far from the boundary.

The expansion transfers source derivatives to the exponential kernels, leaving velocity products and explicit frequency coefficients. For more general finite-energy inputs, it is interpreted by approximation on fixed nonzero frequency bands.  The formulas without boundary terms are used only when the requisite traces vanish.

\subsection{Weighted norms and three operator estimates}
For the tangential argument, define

\begin{equation}
\vartheta_\rho(y)=\chi_\rho(y)^2=\left(\frac\rho{\rho+y}\right)^2,
\qquad V_\rho(s)=\norm{\chi_\rho \,s}_2,
\qquad\rho>0.
\tag{5.11}\label{eq:tangent-weight}
\end{equation}

The density $\vartheta_\rho$ is the square of the multiplier $\chi_\rho$ used in Section 4; the notation is the same as in the introduction. Since $y^2\vartheta_\rho\le\rho^2$, the Hardy inequality gives $V_\rho(s)\le\rho\|s/y\|_2<\infty$. The later absorption argument therefore concerns a quantity known in advance to be finite. The square also matches the factor $\rho/S$ on strips at height $S$.

\result{Lemma 5.3 }
Let $A=(r,s)$ be divergence-free in $\HH$, with zero boundary trace, $A\in L^\infty$, and $\nabla A\in L^2$. Extend $r$ oddly and $s$ evenly. Then

\begin{equation}
\norm{\chi_\rho \,G(A\otimes A)}_2
\le C\bigl(\norm{\nabla A}_2+\norm A_{L^\infty}\bigr)V_\rho(s).
\tag{5.12}\label{eq:weighted-G}
\end{equation}

where $C$ is independent of $\rho$ and $A$.

 The proof is given in Appendix B. The main point is to use $r_x=-s_y$ to control $R_xR_y(r^2)$ by a gradient norm times the norm of the normal component. The other quadratic terms already contain $s$.

\result{Lemma 5.4 }
If $A=(r,s)$ is smooth, compactly supported, divergence-free, and zero on the boundary, then

\begin{equation}
\begin{aligned}
\norm{|D_x|^{-1/2}\Pil T_b(A\otimes A)}_{L_x^2}
\le{}&Cb^{-1}\norm{\nabla A}_2^2\\
&+C\bigl(\norm{\nabla A}_2+\norm A_{L^\infty}\bigr)\norm s_2.
\end{aligned}
\tag{5.13}\label{eq:finite-wall}
\end{equation}

The estimate also holds for $A=\nabla^\perp\phi\in L^\infty(\HH)$ with $\phi\in H_0^2(\HH)$.

\step{Proof.}
Starting from \eqref{eq:wall-reduced}, we estimate the ordinary products, the Oseen square and mixed terms, and the harmonic derivative product separately. Velocity norms are over $\HH$; frequency norms are denoted by $L_k^2$.

\step{Step 1: the five terms to be estimated.}
By \eqref{eq:low-frequency}, on the low-frequency interval,

\[
\frac{\sqrt{|k|}}{|\lambda-|k||}\le Cb^{-1/2},\qquad
\frac{|\lambda|}{|\lambda-|k||}\le C,\qquad
\frac{|k|}{|\lambda-|k||}\le C,\qquad
\frac{|\lambda^2+|k|^2|}{|\lambda|\,|\lambda-|k||}\le C.
\]

Consequently, \eqref{eq:wall-reduced} gives

\[
\begin{aligned}
\norm{|k|^{-1/2}\Pil\widehat{T_b(A\otimes A)}}_{L_k^2}
\le{}&Cb^{-1/2}\norm{\lambda\Tl(r^2)}_{L_k^2}
+C\norm{\lambda|k|^{-1/2}\Tl(rs)}_{L_k^2}\\
&+2\norm{\Pil\frac{\sqrt{|k|}}{\lambda-|k|}\Tk(r_ys)}_{L_k^2}\\
&+C\norm{\sqrt{|k|}\Tl(s^2)}_{L_k^2}
+C\norm{\sqrt{|k|}\Tk(s^2)}_{L_k^2}.
\end{aligned}
\]

We retain the full coefficient in the third term so that $r_ys\in L^1$ can be used directly.

\step{Step 2: traces of ordinary products.}
For $z=|k|$ or $z=\lambda$, Cauchy--Schwarz in the normal variable gives

\[
|\mathcal T_zH(k)|^2
\le\left(\int_0^\infty e^{-2\Re z\,y}\dd y\right)
      \left(\int_0^\infty|\widehat H(k,y)|^2\dd y\right)
=\frac1{2\Re z}\int_0^\infty|\widehat H(k,y)|^2\dd y.
\]

Since $\Re z\ge|k|$, multiplication by $|k|$ and integration in frequency yield

\begin{equation}
\norm{\sqrt{|k|}\,\mathcal T_z(s^2)}_{L_k^2}
\le\frac1{\sqrt2}\norm {(s^2)}_2\le\frac1{\sqrt2}\norm A_{L^\infty}\norm s_2,
\qquad z=|k|,\lambda.
\tag{5.13a}\label{eq:finite-product-trace}
\end{equation}

\step{Step 3: the Oseen square and mixed terms.}
The real part and modulus of $\lambda=\sqrt{|k|^2+ibk}$ satisfy

\[
\Re\lambda
=\left(\frac{\sqrt{k^4+b^2k^2}+k^2}{2}\right)^{1/2}
\ge c\bigl(|k|+\sqrt{b|k|}\bigr),
\qquad
\frac{|\lambda|^2}{|k|}=\sqrt{k^2+b^2}\le |k|+b.
\]

Direct integration over the full frequency axis gives, for every $y>0$,

\[
\begin{aligned}
\int_{\RR}|e^{-y\lambda(k)}|^2\dd k
&\le C\int_0^\infty e^{-c\sqrt{bk}\,y}\dd k
\le \frac{C}{by^2},\\
\int_{\RR}\frac{|\lambda(k)|^2}{|k|}|e^{-y\lambda(k)}|^2\dd k
&\le C\int_0^\infty k e^{-cky}\dd k
   +Cb\int_0^\infty e^{-c\sqrt{bk}\,y}\dd k
\le \frac{C}{y^2}.
\end{aligned}
\]

Fourier inversion and Cauchy--Schwarz imply, for $h\in L^2(\RR)$,

\[
\begin{aligned}
\norm{\mathcal F_x^{-1}(e^{-y\overline\lambda}h)}_{L_x^\infty}
&\le \frac{C}{\sqrt b\,y}\norm h_{L_k^2},\\
\norm{\mathcal F_x^{-1}(\overline\lambda |k|^{-1/2}e^{-y\overline\lambda}h)}_{L_x^\infty}
&\le \frac{C}{y}\norm h_{L_k^2}.
\end{aligned}
\]

The conjugates in these inverse transforms are those appearing in the dual pairings below. For $q\in L^2(\HH)$, the first estimate and the Hardy inequality give

\[
\begin{aligned}
|\langle\Tl(rq),h\rangle_{L_k^2}|
&\le\int_0^\infty\norm{rq(\cdot,y)}_{L_x^1}
   \norm{\mathcal F_x^{-1}(e^{-y\overline\lambda}h)}_{L_x^\infty}\dd y\\
&\le Cb^{-1/2}\norm h_{L_k^2}\int_{\HH}\frac{|r q|}{y}\dd x\dd y\\
&\le Cb^{-1/2}\norm{r/y}_2\norm q_2\norm h_{L_k^2}\\
&\le Cb^{-1/2}\norm{r_y}_2\norm q_2\norm h_{L_k^2}.
\end{aligned}
\]

Taking the supremum over $\|h\|_{L_k^2}=1$, we obtain

\[
\norm{\Tl(rq)}_{L_k^2}
\le Cb^{-1/2}\norm{r_y}_2\norm q_2.
\]

Since $r(x,0)=0$ and the input is compactly supported, integration in the normal variable gives

\[
\lambda\Tl(r^2)
=\mathcal T_\lambda(\partial_y r^2)=2\Tl(rr_y).
\]

Taking $q=r_y$ in the preceding estimate yields

\begin{equation}
\norm{\lambda\Tl(r^2)}_{L_k^2}
\le Cb^{-1/2}\norm{r_y}_2^2
\le Cb^{-1/2}\norm{\nabla A}_2^2.
\tag{5.13b}\label{eq:finite-oseen-square}
\end{equation}

Likewise, the second inverse-transform estimate gives

\[
\begin{aligned}
|\langle\lambda |k|^{-1/2}\Tl(rs),h\rangle_{L_k^2}|
&\le C\norm h_{L_k^2}\int_{\HH}\frac{|rs|}{y}\dd x\dd y\\
&\le C\norm{r/y}_2\norm s_2\norm h_{L_k^2}\\
&\le C\norm{\nabla A}_2\norm s_2\norm h_{L_k^2}.
\end{aligned}
\]

and hence

\begin{equation}
\norm{\lambda |k|^{-1/2}\Tl(rs)}_{L_k^2}
\le C\norm{\nabla A}_2\norm s_2.
\tag{5.13c}\label{eq:finite-oseen-mixed}
\end{equation}

\step{Step 4: the harmonic term.}
By \eqref{eq:low-frequency}, for every $y>0$,

\[
\norm{\Pil\frac{\sqrt{|k|}}{\lambda-|k|}e^{-|k|y}}_{L_k^2}^2
\le\frac Cb\int_{|k|\le b/4}e^{-2|k|y}\dd k\le C.
\]

Minkowski's inequality, the $L^1\to L^\infty$ bound for the horizontal Fourier transform, and Cauchy--Schwarz give

\begin{equation}
\begin{aligned}
\norm{\Pil\frac{\sqrt{|k|}}{\lambda-|k|}\Tk(r_ys)}_{L_k^2}
&\le\int_0^\infty\norm{\Pil\frac{\sqrt{|k|}}{\lambda-|k|}e^{-|k|y}
\widehat{r_ys}(k,y)}_{L_k^2}\dd y\\
&\le C\int_0^\infty\norm{r_y(\cdot,y)s(\cdot,y)}_{L_x^1}\dd y\\
&\le C\norm{r_y}_2\norm s_2.
\end{aligned}
\tag{5.13d}\label{eq:finite-harmonic-square}
\end{equation}

\step{Step 5: combination of the estimates.}
Insert \eqref{eq:finite-product-trace}--\eqref{eq:finite-harmonic-square} into the five-term bound of Step 1 and apply Plancherel's theorem. We obtain

\[
\norm{|D_x|^{-1/2}\Pil T_b(A\otimes A)}_{L_x^2}
\le Cb^{-1}\norm{\nabla A}_2^2
+C\bigl(\norm{\nabla A}_2+\norm A_{L^\infty}\bigr)\norm s_2,
\]

which is \eqref{eq:finite-wall}.

\step{Step 6: extension to $\phi\in H_0^2(\HH)$.}
Let $A=\nabla^\perp\phi\in L^\infty$. First apply a radial cutoff to $\phi$, then translate its zero extension upward and mollify with a nonnegative smooth kernel. This gives $\phi_n\in C_c^\infty(\HH)$ such that $A_n=\nabla^\perp\phi_n$ satisfies

\[
A_n\to A\quad\text{in }H^1(\HH),\qquad
\limsup_{n\to\infty}\norm{A_n}_{L^\infty}\le\norm A_{L^\infty}.
\]

Indeed, the zero extension of $\phi\in H_0^2$ belongs to $H^2(\RR^2)$ and hence to $L^\infty$. The gradient error from the radial cutoff is bounded by $CR^{-1}\|\phi\|_\infty$, while convolution with a nonnegative normalized kernel does not increase the supremum norm. Thus the approximation preserves the required uniform bound, and $A_n\otimes A_n\to A\otimes A$ in $L^2$. On each band $\varepsilon\le|k|\le b/4$, continuity of the source-integral operators passes \eqref{eq:finite-wall} to the limit. Letting $\varepsilon\downarrow0$ and using monotone convergence completes the proof.\qedbox

\result{Lemma 5.5 }
Let $A=(r,s)=\nabla^\perp\phi\in H^1_{\rm loc}(\overline\HH;\RR^2)$ be divergence-free and zero on the boundary, with $\phi(x,0)=0$, $A\in L^\infty(\HH)$, and $\|\nabla A\|_2<\infty$. For $\rho\ge b^{-1}$,

\begin{equation}
\begin{aligned}
\norm{\chi_\rho \,\Py\Pil T_b(A\otimes A)}_2
\le{}&Cb^{-1}\norm{\nabla A}_2^2
+C\bigl(\norm{\nabla A}_2+\norm A_{L^\infty}\bigr)V_\rho(s)\\
&+\frac{C\norm{\nabla A}_2}{\sqrt{b\rho}}V_\rho(s).
\end{aligned}
\tag{5.14}\label{eq:weighted-wall}
\end{equation}

The constant is independent of $\rho$.

\step{Proof.}
We use \eqref{eq:wall-cancel} and split the stream function into a part supported at finite height and the remaining part. Lemma 5.4 handles the first contribution; the second is estimated on dyadic horizontal strips.

\step{Step 1: the sources in the strip decomposition.}
Formula \eqref{eq:wall-cancel} contains three sources, $r_ys$, $rs$, and $s^2$, and five kernels. Each source contains $s$. Since derivatives have already been transferred to the exponential kernels, the products can be decomposed into strips without differentiating the strip cutoffs.

\step{Step 2: a finite-height input.}
Choose a smooth function $0\le\zeta\le1$ equal to one on $[0,1]$ and zero on $[2,\infty)$, and set

\[
\zeta_\rho(y)=\zeta(y/\rho),\qquad
A^{(\rho)}=\nabla^\perp(\zeta_\rho\phi)
=(r^{(\rho)},s^{(\rho)})
=(\zeta_\rho r+\zeta_\rho'\phi,\zeta_\rho s).
\]

Cutting off the stream function preserves incompressibility and the zero boundary trace. Differentiation gives

\[
\begin{aligned}
\partial_xr^{(\rho)}&=\zeta_\rho r_x-\zeta_\rho's,&
\partial_yr^{(\rho)}&=\zeta_\rho r_y+2\zeta_\rho'r+\zeta_\rho''\phi,\\
\partial_xs^{(\rho)}&=\zeta_\rho s_x,&
\partial_ys^{(\rho)}&=\zeta_\rho s_y+\zeta_\rho's.
\end{aligned}
\]

The derivatives of the cutoff are supported in $\rho<y<2\rho$. The Hardy inequalities therefore imply

\[
\begin{aligned}
\norm{\nabla A^{(\rho)}}_2
&\le C\left(\norm{\nabla A}_2
+\rho^{-1}\norm A_{2,\rho<y<2\rho}
+\rho^{-2}\norm\phi_{2,\rho<y<2\rho}\right)\\
&\le C\left(\norm{\nabla A}_2+\norm{A/y}_2+\norm{\phi/y^2}_2\right)
\le C\norm{\nabla A}_2.
\end{aligned}
\]

Also $|\phi(x,y)|\le y\|A\|_\infty$, so

\[
\norm{A^{(\rho)}}_{L^\infty}
\le\norm A_{L^\infty}+\frac C\rho\norm\phi_{\infty,\rho<y<2\rho}
\le C\norm A_{L^\infty}.
\]

Since $\chi_\rho \ge1/3$ for $0<y<2\rho$,

\[
\norm{s^{(\rho)}}_2\le3V_\rho(s).
\]

Finite-height support and $\phi/y^2,A/y,\nabla A\in L^2$ imply $\zeta_\rho\phi\in H_0^2(\HH)$. The preceding supremum bound allows us to apply the extension of Lemma 5.4. Using $\vartheta_\rho\le1$ and the Poisson identity

\[
\norm{\Py q}_{L^2(\HH)}^2
=\int_\RR\int_0^\infty e^{-2|k|y}|\widehat q(k)|^2\dd y\dd k
=\frac12\norm{|D_x|^{-1/2}q}_{L_x^2}^2,
\]

we obtain

\[
\begin{aligned}
&\norm{\chi_\rho \,\Py\Pil
 T_b(A^{(\rho)}\otimes A^{(\rho)})}_2\\
&\qquad\le Cb^{-1}\norm{\nabla A^{(\rho)}}_2^2
+C\bigl(\norm{\nabla A^{(\rho)}}_2+\norm{A^{(\rho)}}_{L^\infty}\bigr)
 \norm{s^{(\rho)}}_2\\
&\qquad\le Cb^{-1}\norm{\nabla A}_2^2
+C\bigl(\norm{\nabla A}_2+\norm A_{L^\infty}\bigr)V_\rho(s).
\end{aligned}
\]

It remains to estimate the difference between the boundary responses of $A$ and $A^{(\rho)}$.

\step{Step 3: kernels for sources at large height.}
Suppose $S<t<2S$, with $S\ge\rho\ge b^{-1}$. All pointwise bounds in the table are on $0<|k|\le b/4$; the final column bounds $\|\Pil K(\cdot,t)\|_{L_k^2}$. The estimates follow by direct integration from \eqref{eq:low-frequency} and $\Re\lambda\ge c\sqrt{b|k|}$. Constants are independent of $S,t,\rho,b$.
\begin{center}
\renewcommand{\arraystretch}{2.1}
\setlength{\tabcolsep}{4pt}
\small
\begin{tabular}{cccc}
\toprule
Source & Kernel $K(k,t)$ & Pointwise bound & $L_k^2$ bound\\
\midrule
$r_ys$ & $\displaystyle\frac{2\lambda e^{-\lambda t}}{\lambda-|k|}$
& $Ce^{-cS\sqrt{b|k|}}$ & $\displaystyle\frac C{\sqrt b\,S}$\\
$r_ys$ & $\displaystyle\frac{-2|k|e^{-|k|t}}{\lambda-|k|}$
& $Cb^{-1/2}|k|^{1/2}e^{-S|k|}$ & $\displaystyle\frac C{\sqrt b\,S}$\\
$rs$ & $-(\lambda+|k|)e^{-\lambda t}$
& $C\sqrt b\,|k|^{1/2}e^{-cS\sqrt{b|k|}}$ & $\displaystyle\frac C{\sqrt b\,S^2}$\\
$s^2$ & $\displaystyle\frac{ik|k|e^{-|k|t}}{\lambda-|k|}$
& $Cb^{-1/2}|k|^{3/2}e^{-S|k|}$ & $\displaystyle\frac C{\sqrt b\,S^2}$\\
$s^2$ & $\displaystyle\frac{-ik\lambda e^{-\lambda t}}{\lambda-|k|}$
& $C|k|e^{-cS\sqrt{b|k|}}$ & $\displaystyle\frac C{b^{3/2}S^3}\le\frac C{\sqrt b\,S^2}$\\
\bottomrule
\end{tabular}
\end{center}
The last row uses $bS\ge1$. Thus the two kernels acting on $r_ys$ have norm at most $C/(\sqrt b\,S)$, and the other three have norm at most $C/(\sqrt b\,S^2)$.

\step{Step 4: one strip of the boundary field.}
Set

\[
\Sigma_S=\{(x,y):S<y<2S\},\qquad
V_S=\norm s_{L^2(\Sigma_S)},\qquad
d_S^2=\int_{\Sigma_S}\left(|r_y|^2+\frac{|A|^2}{y^2}
+\frac{\phi^2}{y^4}\right).
\]

By definition,

\[
\norm{r_y}_{2,\Sigma_S}\le d_S,\qquad
\norm A_{2,\Sigma_S}\le2Sd_S,\qquad
\norm\phi_{2,\Sigma_S}\le4S^2d_S.
\]

Cauchy--Schwarz gives

\[
\begin{aligned}
\norm{r_ys}_{1,\Sigma_S}&\le d_SV_S,\\
\norm{rs}_{1,\Sigma_S}+\norm{s^2}_{1,\Sigma_S}&\le CSd_SV_S,\\
\norm{\phi s}_{1,\Sigma_S}&\le CS^2d_SV_S.
\end{aligned}
\]

For a source $H$ supported in $\Sigma_S$, write its boundary contribution as

\[
\widehat q(k)=\Pil\int_S^{2S}K(k,t)\widehat H(k,t)\dd t.
\]

The horizontal Fourier $L^1\to L^\infty$ bound and Minkowski's inequality yield

\[
\begin{aligned}
\norm{\widehat q}_{L_k^2}
&\le\int_S^{2S}\norm{\Pil K(\cdot,t)\widehat H(\cdot,t)}_{L_k^2}\dd t\\
&\le C\sup_{S<t<2S}\norm{\Pil K(\cdot,t)}_{L_k^2}
       \int_S^{2S}\norm{H(\cdot,t)}_{L_x^1}\dd t\\
&=C\sup_{S<t<2S}\norm{\Pil K(\cdot,t)}_{L_k^2}\norm H_{1,\Sigma_S}.
\end{aligned}
\]

Since $\int_0^\infty\vartheta_\rho(y)\dd y=\rho$, Plancherel gives

\[
\begin{aligned}
\norm{\chi_\rho \,\Py q}_2^2
&=\int_\RR|\widehat q(k)|^2
  \left(\int_0^\infty\vartheta_\rho(y)e^{-2|k|y}\dd y\right)\dd k\\
&\le\rho\norm{\widehat q}_{L_k^2}^2.
\end{aligned}
\]

Thus the weighted Poisson field generated by this source is bounded by

\[
C\sqrt\rho\sup_{S<t<2S}\norm{\Pil K(\cdot,t)}_{L_k^2}\norm H_{1,\Sigma_S}.
\]

Apply the kernel bound $C/(\sqrt b\,S)$ to $r_ys$ and $C/(\sqrt b\,S^2)$ to $rs$ and $s^2$. With the product estimates above, the total contribution of a single strip is at most

\[
\frac{C\sqrt\rho}{\sqrt b\,S}d_SV_S
=\frac C{\sqrt{b\rho}}d_S\frac\rho S V_S.
\]

\step{Step 5: cutoff errors and summation.}
Apply \eqref{eq:wall-cancel} to both $A$ and $A^{(\rho)}$. The differences at the three source positions are

\[
\begin{aligned}
r_ys-\partial_yr^{(\rho)}s^{(\rho)}
&=(1-\zeta_\rho^2)r_ys
  -2\zeta_\rho\zeta_\rho'rs-\zeta_\rho\zeta_\rho''\phi s,\\
rs-r^{(\rho)}s^{(\rho)}
&=(1-\zeta_\rho^2)rs-\zeta_\rho\zeta_\rho'\phi s,\\
s^2-(s^{(\rho)})^2&=(1-\zeta_\rho^2)s^2.
\end{aligned}
\]

All these sources are supported in $y\ge\rho$. On the first strip, $S=\rho$, we have $|\zeta_\rho'|\le C/\rho$ and $|\zeta_\rho''|\le C/\rho^2$, and hence

\[
\begin{aligned}
\norm{2\zeta_\rho\zeta_\rho'rs}_{1,\Sigma_\rho}
&\le C\rho^{-1}(\rho d_S V_S)
\le Cd_S V_S,\\
\norm{\zeta_\rho\zeta_\rho''\phi s}_{1,\Sigma_\rho}
&\le C\rho^{-2}(\rho^2d_S V_S)
\le Cd_S V_S,\\
\norm{\zeta_\rho\zeta_\rho'\phi s}_{1,\Sigma_\rho}
&\le C\rho^{-1}(\rho^2d_S V_S)
\le C\rho d_S V_S.
\end{aligned}
\]

On the remaining strips $\Sigma_{2^j\rho}$, $j\ge1$, the cutoff derivatives vanish. Thus, for every $S=2^j\rho$, the first difference source has $L^1$ norm bounded by $Cd_SV_S$, and each of the other two by $CSd_SV_S$. Step 4 therefore applies to every difference source.

The strips are disjoint. By the Hardy inequalities,

\[
\sum_{S=2^j\rho,\ j\ge0}d_S^2
\le\norm{r_y}_2^2+\norm{A/y}_2^2+\norm{\phi/y^2}_2^2
\le C\norm{\nabla A}_2^2.
\]

Also, for $S\ge\rho$ and $S<y<2S$,

\[
\frac\rho S\le3\frac\rho{\rho+y}=3\chi_\rho(y).
\]

so that

\[
\sum_{S=2^j\rho,\ j\ge0}\left(\frac\rho S V_S\right)^2
\le9\int_{y>\rho}\vartheta_\rho(y)|s(x,y)|^2\dd x\dd y
\le9V_\rho(s)^2.
\]

The triangle inequality and Cauchy--Schwarz now give

\[
\begin{aligned}
&\norm{\chi_\rho \,\Py\Pil
 \{T_b(A\otimes A)-T_b(A^{(\rho)}\otimes A^{(\rho)})\}}_2\\
&\quad\le\frac C{\sqrt{b\rho}}\sum_{S=2^j\rho}d_S\frac\rho S V_S\\
&\quad\le\frac C{\sqrt{b\rho}}
 \left(\sum_Sd_S^2\right)^{1/2}
 \left(\sum_S\left(\frac\rho S V_S\right)^2\right)^{1/2}\\
&\quad\le\frac{C\norm{\nabla A}_2}{\sqrt{b\rho}}V_\rho(s).
\end{aligned}
\]

This also proves absolute convergence of the strip contributions in weighted $L^2$. Adding the finite-height estimate from Step 2 yields \eqref{eq:weighted-wall}.\qedbox

% section6
\section{Tangential backgrounds: absorption and rigidity}
We now combine the two identities and three estimates from Section 5. The order of the argument is important: first bound the finite quantity $V_\rho(v)$ uniformly in $\rho$, with the compact part fixed, and only then let $\rho\to\infty$.

\result{Theorem 6.1 }
Under the assumptions of Theorem 1.1, if $B=(b,0)$ with $b\neq0$, then $v\in L^2(\HH)$.

\step{Proof.}
By Lemma 5.1, the solution satisfies

\[
bv=\omega-h_f-G(W\otimes W),\qquad h_f=\Py f,\quad f=\omega(\cdot,0).
\]

We estimate this identity in weighted $L^2$. The quadratic source will be split into a fixed compact part and a small divergence-free tail. The tail contributes at most a small multiple of $V_\rho(v)$, while the compact contribution is bounded independently of $\rho$.

\step{Step 1: a small divergence-free tail.}
Write $\varrho=(x^2+y^2)^{1/2}$. Lemma 2.1 gives $W/\varrho,\psi/\varrho^2\in L^2$. Moreover, $W\in C_0$ implies uniform decay of $\psi/\varrho$. Indeed, given $\varepsilon>0$, choose a fixed disk outside which $|u|\le\varepsilon$.  Thus the formula $\psi(x,y)=\int_0^yu(x,t)\dd t$ gives

\[
|\psi(x,y)|\le C_\varepsilon+\varepsilon y,
\qquad\limsup_{\varrho\to\infty}|\psi|/\varrho\le\varepsilon.
\]

Choose a radial cutoff $\eta_R$ equal to one on $B_R$ and zero outside $B_{2R}$, with $|D^j\eta_R|\le C_jR^{-j}$. Set

\[
\psi_c=\eta_R\psi,\qquad \phi=(1-\eta_R)\psi,
\qquad W_c=\nabla^\perp\psi_c,\qquad A=\nabla^\perp\phi=(r,s).
\]

Then $W=W_c+A$, and both fields are divergence-free and zero on the boundary. The product rule gives

\[
\begin{aligned}
\norm{\nabla A}_2
&\le C\left(\norm{\nabla W}_{2,\varrho>R}+\norm{W/\varrho}_{2,\varrho>R}
+\norm{\psi/\varrho^2}_{2,\varrho>R}\right),\\
\norm A_{L^\infty}
&\le C\left(\norm W_{\infty,\varrho>R}+\norm{\psi/\varrho}_{\infty,\varrho>R}\right).
\end{aligned}
\]

We can therefore fix $R$ sufficiently large that

\begin{equation}
\norm{\nabla A}_2+\norm A_{L^\infty}\le\varepsilon.
\tag{6.1}
\end{equation}

From now on $R$ is fixed, and every constant associated with the compact part is independent of $\rho$. Define

\[
K=W\otimes W-A\otimes A
=W_c\otimes W_c+W_c\otimes A+A\otimes W_c.
\]

\step{Step 2: the bulk contribution.}
Write $V_\rho=V_\rho(v)$. Since $s=v-(W_c)_2$,

\[
V_\rho(s)\le V_\rho+\norm{(W_c)_2}_2.
\]

Lemma 5.3 and (6.1) imply

\[
\begin{aligned}
\norm{\chi_\rho G(A\otimes A)}_2
&\le C(\norm{\nabla A}_2+\norm A_\infty)V_\rho(s)\\
&\le C\varepsilon V_\rho+C\varepsilon\norm{(W_c)_2}_2.
\end{aligned}
\]

The parity extension of $K$ belongs to $L^2(\RR^2)$. By Riesz-transform boundedness and $0<\vartheta_\rho\le1$,

\[
\norm{\chi_\rho G(K)}_2\le\norm{G(K)}_2\le C\norm K_2.
\]

Using linearity of $G$ and including the fixed compact terms in $C_R$, we obtain

\begin{equation}
\norm{\chi_\rho \,G(W\otimes W)}_2
\le C\varepsilon V_\rho+C_R.
\tag{6.2}\label{eq:core-bulk}
\end{equation}

\step{Step 3: the boundary contribution.}
Lemma 5.2 and linearity of $T_b$ give at low frequencies

\[
\Pil f=\Pil T_b(K)+\Pil T_b(A\otimes A).
\]

We estimate the compact and tail contributions at low frequencies, then handle the high-frequency part of the actual boundary data. Throughout this step, $\rho\ge b^{-1}$.

\textit{(i) The compact contribution at low frequencies.} For each entry $K_{ij}$ and for $z=|k|,\lambda$, the inequality $\Re z>0$ gives

\[
\begin{aligned}
|\mathcal T_zK_{ij}(k)|
&=\left|\int_0^\infty e^{-zy}\widehat K_{ij}(k,y)\dd y\right|\\
&\le C\int_0^\infty\int_\RR|K_{ij}(x,y)|\dd x\dd y
\le C\norm K_1.
\end{aligned}
\]

The tensor $K$ is a difference of quadratic tensors formed from zero-boundary fields, so \eqref{eq:wall-expansion} applies to it. On $0<|k|\le b/4$, we have $|\lambda|\asymp\sqrt{b|k|}$ and $|\lambda-|k||\asymp\sqrt{b|k|}$. Hence the coefficients in that expansion satisfy

\[
\frac{|\lambda^2+|k|^2|+2|k|^2+|k\lambda|+|k|^2}
{\sqrt{|k|}\,|\lambda-|k||}
\le C\sqrt b.
\]

Combining the last two bounds gives

\[
|k|^{-1/2}|\widehat{T_b(K)}(k)|\le C\sqrt b\norm K_1.
\]

The low-frequency interval has total length $b/2$, and therefore

\[
\int_{0<|k|\le b/4}\frac{|\widehat{T_b(K)}(k)|^2}{|k|}\dd k
\le Cb^2\norm K_1^2.
\]

Since $\widehat{\Py g}(k)=e^{-|k|y}\widehat g(k)$, Plancherel's theorem yields

\begin{equation}
\begin{aligned}
\norm{\chi_\rho \Py\Pil T_b(K)}_2
&\le\norm{\Py\Pil T_b(K)}_2\\
&=\left(\int_{0<|k|\le b/4}
\frac{|\widehat{T_b(K)}(k)|^2}{2|k|}\dd k\right)^{1/2}
\le Cb\norm K_1.
\end{aligned}
\tag{6.3}\label{eq:core-wall}
\end{equation}

The bound is independent of $\rho$, since $\vartheta_\rho\le1$.

\textit{(ii) The tail at low frequencies.} Apply \eqref{eq:weighted-wall} from Lemma 5.5 to $A=(r,s)$:

\[
\begin{aligned}
\norm{\chi_\rho \Py\Pil T_b(A\otimes A)}_2
\le{}&Cb^{-1}\norm{\nabla A}_2^2\\
&+C\left(\norm{\nabla A}_2+\norm A_\infty
+\frac{\norm{\nabla A}_2}{\sqrt{b\rho}}\right)V_\rho(s).
\end{aligned}
\]

The condition $\rho\ge b^{-1}$ gives $(b\rho)^{-1/2}\le1$. Using (6.1) and $V_\rho(s)\le V_\rho+\|(W_c)_2\|_2$, we conclude that

\[
\norm{\chi_\rho \Py\Pil T_b(A\otimes A)}_2
\le Cb^{-1}\varepsilon^2+C\varepsilon V_\rho
+C\varepsilon\norm{(W_c)_2}_2.
\]

Together with \eqref{eq:core-wall}, this gives the complete low-frequency bound

\[
\norm{\chi_\rho \Py\Pil f}_2
\le C\varepsilon V_\rho
+C\left(b\norm K_1+b^{-1}\varepsilon^2
+\varepsilon\norm{(W_c)_2}_2\right).
\]

\textit{(iii) High frequencies.} For $f\in L^2(\RR)$ and  $\vartheta_\rho\le1$, using Plancherel's theorem and Tonelli's theorem, we obtain

\[
\begin{aligned}
\norm{\chi_\rho \Py(1-\Pil)f}_2^2
&\le\int_{|k|>b/4}|\widehat f(k)|^2
\left(\int_0^\infty e^{-2|k|y}\dd y\right)\dd k\\
&=\int_{|k|>b/4}\frac{|\widehat f(k)|^2}{2|k|}\dd k
\le\frac2b\norm f_2^2.
\end{aligned}
\]

Combining the two frequency ranges in $h_f=\Py\Pil f+\Py(1-\Pil)f$ gives

\begin{equation}
\norm{\chi_\rho \,h_f}_2
\le C\varepsilon V_\rho+C_R.
\tag{6.4}
\end{equation}

The constant $C_R$ includes $b\|K\|_1$, $b^{-1}\varepsilon^2$, $\varepsilon\|(W_c)_2\|_2$, and $b^{-1/2}\|f\|_2$.

\step{Step 4: absorption before removal of the weight.}
Taking the weighted $L^2$ norm in \eqref{eq:tangent-identity} and using (6.2) and (6.4), we obtain

\[
bV_\rho\le\norm\omega_2+C_R+C\varepsilon V_\rho.
\]

Choose $\varepsilon$ so that $C\varepsilon\le b/2$, and then fix the corresponding radius $R$ from Step 1. The Hardy inequality already guarantees $V_\rho<\infty$, so absorption gives

\[
V_\rho\le\frac Cb(\norm\omega_2+C_R),\qquad\rho\ge b^{-1}.
\]

The right-hand side is independent of $\rho$. Since $\vartheta_\rho\uparrow1$, monotone convergence as $\rho\to\infty$ proves $v\in L^2(\HH)$.\qedbox

\subsection*{A truncated stream-function test}
\result{Theorem 6.2 (Rigidity for a tangential background).}
Under the assumptions of Theorem 1.1, if $B=(b,0)$ with $b\ne0$, then $U\equiv B$ in $\HH$ and $P$ is constant.

\step{Proof. Step 1: a test preserving the nonlinear cancellation.}
A horizontal reflection allows us to assume $b>0$. By Theorem 6.1, $v\in L^2(\HH)$. Choose $\chi_S(y)=\chi(y/S)$, where $\chi=1$ on $[0,1]$ and $\chi=0$ on $[2,\infty)$, and set $\varphi_S=\chi_S\psi$. Finite-height support and the bounds $W/y,\psi/y^2,D^2\psi\in L^2$ imply $\varphi_S\in H_0^2(\HH)$. The following distributional vorticity identity holds:

\[
-\int\omega\Delta\varphi-\int\omega(B+W)\cdot\nabla\varphi=0.
\]

 Since $W\cdot\nabla\psi=0$,

\[
\Delta\varphi_S=-\chi_S\omega+2\chi_S'u+\chi_S''\psi,
\qquad
(B+W)\cdot\nabla\varphi_S=-b\chi_Sv+\chi_S'\psi v.
\]

we obtain

\begin{equation}
\int\chi_S\omega^2
=-b\int\chi_S\omega v+\int\chi_S'\omega\psi v
+2\int\chi_S'\omega u+\int\chi_S''\omega\psi.
\tag{6.5}
\end{equation}

Incompressibility gives $\omega v=\frac12\partial_x(v^2-u^2)-\partial_y(uv)$, hence $\int\chi_S\omega v=\int\chi_S'uv$. These integrations can first be performed with a horizontal cutoff and then justified using square integrability on finite-height strips and the zero boundary trace.

\step{Step 2: passage to infinite height.}
Let $\Sigma_S=\{(x,y)\in\HH:S<y<2S\}$. After the preceding substitution, the four terms on the right of (6.5) satisfy

\[
\begin{aligned}
\left|b\int_{\HH}\chi_S'uv\right|
&\le C|b|\norm{u/y}_{2,\Sigma_S}\norm v_{2,\Sigma_S},\\
\left|\int_{\HH}\chi_S'\omega\psi v\right|
&\le C\norm u_\infty\norm\omega_{2,\Sigma_S}\norm v_{2,\Sigma_S},\\
\left|2\int_{\HH}\chi_S'\omega u\right|
&\le C\norm\omega_{2,\Sigma_S}\norm{u/y}_{2,\Sigma_S},\\
\left|\int_{\HH}\chi_S''\omega\psi\right|
&\le C\norm\omega_{2,\Sigma_S}\norm{\psi/y^2}_{2,\Sigma_S}.
\end{aligned}
\]

We used $|\chi_S'|\le C/y$, $|\chi_S''|\le C/y^2$ on $\Sigma_S$, and $|\psi(x,y)|\le\|u\|_\infty y$. Since $v,\omega,u/y,\psi/y^2\in L^2(\HH)$, all four terms tend to zero. The dominated convergence theorem on the left gives $\|\omega\|_2^2=E=0$. Thus $\nabla W=0$, and the zero trace implies $W=0$.  Then $U\equiv B$ in $\HH$ and $P$ is constant.\qedbox

% appendices
\appendix
\section{An anisotropic endpoint estimate for Riesz transforms}
All functions and norms in this appendix are on $\RR^2$. Write $\mathcal R=R_xR_y$, with Fourier symbol $-k\eta/(k^2+\eta^2)$. Its value at the origin is immaterial for the $L^2$ operator. For the classical multiplier and singular-integral descriptions of Riesz transforms, see \cite{Stein1970}. We first prove a product estimate and then use incompressibility to obtain the quadratic estimate needed in Appendix B.

\result{Lemma A.1 (A product estimate and a partial Fourier kernel).}
If $f,g,f_x,g_x\in L^2(\RR^2)$, then

\begin{equation}
\norm{\mathcal R(fg)}_2
\le C\bigl(\norm f_2\norm{g_x}_2+\norm{f_x}_2\norm g_2\bigr).
\tag{A.1}
\end{equation}

\step{Proof. Step 1: a one-dimensional product estimate.}
In this step all norms are on $\RR$, and $|D|^{1/2}$ is the multiplier $|k|^{1/2}$. For real-valued $h\in H^1(\RR)$, Plancherel's theorem and the Cauchy--Schwarz inequality give

\[
\begin{aligned}
\norm{|D|^{1/2}(h^2)}_2^2
&=\int_{\RR}|k|\,|\widehat{h^2}(k)|^2\dd k
\le\norm{h^2}_2\norm{(h^2)'}_2\\
&\le2\norm h_\infty^2\norm h_2\norm{h'}_2
\le C\norm h_2^2\norm{h'}_2^2.
\end{aligned}
\]

The last step uses $\|h\|_\infty^2\le2\|h\|_2\|h'\|_2$, which follows from $h(x)^2=2\int_{-\infty}^xh(t)h'(t)\dd t$. For real-valued $f,g\in H^1(\RR)$ and $t>0$, the identity

\[
4fg=(tf+t^{-1}g)^2-(tf-t^{-1}g)^2,
\]

implies

\[
\begin{aligned}
\norm{|D|^{1/2}(fg)}_2\le C\bigl(&t^2\norm f_2\norm{f'}_2
+t^{-2}\norm g_2\norm{g'}_2\\
&+\norm f_2\norm{g'}_2+\norm{f'}_2\norm g_2\bigr).
\end{aligned}
\]

The infimum of the first two terms over $t>0$ is $2(\|f\|_2\|f'\|_2\|g\|_2\|g'\|_2)^{1/2}$, which is bounded by the sum of the two cross terms. Therefore

\begin{equation}
\norm{|D|^{1/2}(fg)}_2
\le C\bigl(\norm f_2\norm{g'}_2+\norm{f'}_2\norm g_2\bigr).
\tag{A.2}
\end{equation}
Complex-valued functions follow by separating real and imaginary parts.

\step{Step 2: the convolution kernel in the normal variable.}
Initially suppose $f,g\in C_c^\infty(\RR^2)$ and set $F=fg$. Write $\widehat F=\mathcal F_xF$ and $\widetilde F=\mathcal F_y\widehat F$. By the multiplier definition,

\[
\widetilde{\mathcal RF}(k,\eta)
=\mathcal F_y\bigl[\widehat{\mathcal RF}(k,\cdot)\bigr](\eta)
=-\frac{k\eta}{k^2+\eta^2}
  \mathcal F_y\bigl[\widehat F(k,\cdot)\bigr](\eta).
\]

For fixed $k\ne0$, inverse transformation in the normal frequency gives

\[
\begin{aligned}
\widehat{\mathcal RF}(k,y)
&=\mathcal F_y^{-1}\left[
  -\frac{k\eta}{k^2+\eta^2}
  \mathcal F_y\bigl[\widehat F(k,\cdot)\bigr](\eta)\right](y)\\
&=\bigl(K_k*_y\widehat F(k,\cdot)\bigr)(y)
=\int_{\RR}K_k(y-z)\widehat F(k,z)\dd z,
\end{aligned}
\]

where

\[
K_k(t)=-\frac{ik}{2}\operatorname{sgn}(t)e^{-|k||t|}.
\]

\step{Step 3: the two-dimensional product estimate.}
Since $\|K_k\|_{L_t^2}^2=|k|/4$, Minkowski's integral inequality and the horizontal Plancherel theorem (see \cite[Sections~1.2 and~2.2]{Grafakos2014}) give

\[
\begin{aligned}
\norm{\mathcal RF}_2
&\le\int_{\RR}\left(\int_{\RR}\norm{K_k}_{L_y^2}^2
             |\widehat F(k,z)|^2\dd k\right)^{1/2}\dd z\\
&=\frac12\int_{\RR}\norm{|D_x|^{1/2}F(\cdot,z)}_{L_x^2}\dd z.
\end{aligned}
\]

Apply (A.2) on each slice and then Cauchy--Schwarz in $z$:

\[
\begin{aligned}
\norm{\mathcal R(fg)}_2
&\le C\int_{\RR}\bigl(
\norm{f(\cdot,z)}_{L_x^2}\norm{g_x(\cdot,z)}_{L_x^2}
+\norm{f_x(\cdot,z)}_{L_x^2}\norm{g(\cdot,z)}_{L_x^2}\bigr)\dd z\\
&\le C\bigl(\norm f_2\norm{g_x}_2+\norm{f_x}_2\norm g_2\bigr).
\end{aligned}
\]

\step{Step 4: approximation.}
Truncation and mollification give $f_m,g_m\in C_c^\infty(\RR^2)$ with $(f_m,\partial_xf_m)\to(f,f_x)$ and $(g_m,\partial_xg_m)\to(g,g_x)$ in $L^2\times L^2$. Applying the estimate to $f_mg_m-f_ng_n=(f_m-f_n)g_m+f_n(g_m-g_n)$ shows that $\mathcal R(f_mg_m)$ is Cauchy in $L^2$. On the other hand, $f_mg_m\to fg$ in $L^1$, so their Fourier transforms converge uniformly. Multiplication by the bounded symbol $-k\eta/(k^2+\eta^2)$ preserves convergence in distributions. The $L^2$ limit is therefore $\mathcal R(fg)$, and (A.1) follows.\qedbox

\result{Remark A.1 }
The same argument gives, for $f,g,f_x,g_x\in L^2(\RR^2)$,

\[
\norm{R_x^2(fg)}_2
\le C\bigl(\norm f_2\norm{g_x}_2+\norm{f_x}_2\norm g_2\bigr).
\]

Interchanging $x$ and $y$ gives, when $f,g,f_y,g_y\in L^2(\RR^2)$,

\[
\norm{R_y^2(fg)}_2
\le C\bigl(\norm f_2\norm{g_y}_2+\norm{f_y}_2\norm g_2\bigr).
\]

\result{Lemma A.2 }
If $(r,s)\in H^1(\RR^2;\RR^2)$ and $r_x+s_y=0$ , then

\begin{equation}
\norm{R_xR_y(r^2)}_2\le C\norm{\nabla r}_2\norm s_2.
\tag{A.3}
\end{equation}

\step{Proof.}
We use Bony's decomposition \cite{Bony1981} in the normal frequency, separating comparable and widely separated input frequencies; see also \cite{BCD2011}. Products of comparable frequencies are handled by Lemma A.1, without restricting their output frequency. For separated frequencies, the product remains at the larger normal frequency; the multiplier $\mathcal R$ then supplies its reciprocal. Incompressibility transfers the horizontal derivative of $r$ to a normal derivative of $s$, producing a summable factor in the frequency gap.

\step{Step 1: a homogeneous dyadic decomposition.}
We use a one-dimensional dyadic partition of unity in the normal variable, as in \cite{BCD2011}. Choose a real even function $\chi\in C_c^\infty(\RR)$, nonincreasing in $|\eta|$, with $0\le\chi\le1$, equal to one for $|\eta|\le1$ and zero for $|\eta|\ge2$. Set

\[
\varphi(\eta)=\chi(\eta)-\chi(2\eta),\qquad
\varphi_j(\eta)=\varphi(2^{-j}\eta),\quad j\in\mathbb Z.
\]

Then $\varphi$ is real, even, and nonnegative, and

\[
\operatorname{supp}\varphi_j
\subset\{2^{j-1}\le|\eta|\le2^{j+1}\},\qquad
\sum_{j\in\mathbb Z}\varphi_j(\eta)=1\quad(\eta\ne0).
\]

Let $\mathcal F_y$ denote the normal Fourier transform and define

\[
\Delta_j^yh=\mathcal F_y^{-1}\bigl(\varphi_j\mathcal F_yh\bigr),
\qquad r_j=\Delta_j^yr,\quad s_j=\Delta_j^ys.
\]

They commute with both derivatives; hence $r_{j,x}+s_{j,y}=0$.  By Plancherel's theorem and the dominated convergence theorem, for $h\in H^1(\RR^2)$,

\[
\sum_j\norm{\Delta_j^yh}_2^2\le C\norm h_2^2,
\qquad
\sum_{j=-N}^N\Delta_j^yh\longrightarrow h\quad\text{in }H^1(\RR^2).
\]

\step{Step 2: estimates on a single band.}
Write

\[
A_j=\norm{\nabla r_j}_2,\qquad
V_j=\norm{s_j}_2,\qquad B_j=(A_jV_j)^{1/2}.
\]

On the $j$th band, $2^{j-1}\le|\eta|\le2^{j+1}$. Plancherel's theorem and incompressibility give

\begin{equation}
\begin{aligned}
\norm{r_j}_2&\le2^{1-j}\norm{r_{j,y}}_2\le C2^{-j}A_j,\\
\norm{r_{j,x}}_2&=\norm{s_{j,y}}_2\le C2^jV_j,\\
V_j&\le2^{1-j}\norm{s_{j,y}}_2
=2^{1-j}\norm{r_{j,x}}_2\le C2^{-j}A_j.
\end{aligned}
\tag{A.4}
\end{equation}

In particular, $V_j^2\le C2^{-j}A_jV_j$, so $V_j\le C2^{-j/2}B_j$.

We also need a supremum bound for $r_j$. For almost every $x$, the one-dimensional Bernstein inequality in the normal variable (see \cite{BCD2011}), or directly Fourier inversion and Cauchy--Schwarz, gives

\[
\norm{r_j(x,\cdot)}_{L_y^\infty}
\le C2^{j/2}\norm{r_j(x,\cdot)}_{L_y^2}.
\]

Let $q(x)=\|r_j(x,\cdot)\|_{L_y^2}^2$. The differentiation rule for $L_y^2$-valued $H_x^1$ functions yields

\[
|q'(x)|\le2\norm{r_j(x,\cdot)}_{L_y^2}
                  \norm{r_{j,x}(x,\cdot)}_{L_y^2},\qquad
\sup_xq(x)\le\int_{\RR}|q'(x)|\dd x
\le2\norm{r_j}_2\norm{r_{j,x}}_2.
\]

Since $q\in W^{1,1}(\RR)\cap L^1(\RR)$, its limits at infinity vanish. Combining this with (A.4),

\begin{equation}
\norm{r_j}_\infty
\le C2^{j/2}(\norm{r_j}_2\norm{r_{j,x}}_2)^{1/2}
\le C2^{j/2}B_j.
\tag{A.5}
\end{equation}

Finite overlap of the frequency bands and Cauchy--Schwarz imply

\begin{equation}
\begin{aligned}
\sum_jA_j^2&\le C\norm{\nabla r}_2^2,\qquad
\sum_jV_j^2\le C\norm s_2^2,\\
\sum_jB_j^2&=\sum_jA_jV_j
\le\left(\sum_jA_j^2\right)^{1/2}
    \left(\sum_jV_j^2\right)^{1/2}
\le C\norm{\nabla r}_2\norm s_2.
\end{aligned}
\tag{A.6}
\end{equation}

\step{Step 3: comparable frequencies, $|j-\ell|\le3$.}
Lemma A.1 and (A.4) give

\[
\begin{aligned}
\norm{\mathcal R(r_jr_\ell)}_2
&\le C\bigl(\norm{r_j}_2\norm{r_{\ell,x}}_2
            +\norm{r_{j,x}}_2\norm{r_\ell}_2\bigr)\\
&\le C\bigl(2^{\ell-j}A_jV_\ell+2^{j-\ell}V_jA_\ell\bigr)
\le C(A_jV_\ell+V_jA_\ell).
\end{aligned}
\]

For a fixed gap $d=\ell-j\in\{-3,\ldots,3\}$,

\[
\sum_jA_jV_{j+d}
\le\left(\sum_jA_j^2\right)^{1/2}
    \left(\sum_jV_{j+d}^2\right)^{1/2},
\]

The other term is treated in the same way. Summing over these seven gaps and using (A.6), we obtain

\begin{equation}
\sum_{|j-\ell|\le3}\norm{\mathcal R(r_jr_\ell)}_2
\le C\norm{\nabla r}_2\norm s_2.
\tag{A.7}
\end{equation}

No lower bound was placed on the output frequency. Thus the estimate includes low frequencies produced by cancellation of two high-frequency inputs.

\step{Step 4: separated frequencies, $\ell\ge j+4$.}
The normal Fourier support of the product lies in the sum of the two input supports. Consequently,

\[
2^{\ell-2}\le|\eta|\le2^{\ell+2}.
\]

Choose a real even function $\vartheta\in C_c^\infty(\RR)$ equal to one on $1/4\le|\eta|\le4$ and supported in $1/8\le|\eta|\le8$. Define the whole-plane multiplier

\[
\widetilde{\mathcal S_\ell h}(k,\eta)
=\frac{i\eta}{k^2+\eta^2}\vartheta(2^{-\ell}\eta)\widetilde h(k,\eta).
\]

On its support, $|\eta|\ge2^{\ell-3}$. Hence

\[
\norm{\mathcal S_\ell}_{L^2\to L^2}\le C2^{-\ell},
\qquad
\mathcal R(r_jr_\ell)
=\mathcal S_\ell\partial_x(r_jr_\ell).
\]

The product rule and (A.4)--(A.5) give

\[
\begin{aligned}
\norm{\mathcal R(r_jr_\ell)}_2
&\le C2^{-\ell}\bigl(
\norm{r_j}_\infty\norm{r_{\ell,x}}_2
+\norm{r_\ell}_\infty\norm{r_{j,x}}_2\bigr)\\
&\le C\bigl(2^{j/2}B_jV_\ell+2^{j-\ell/2}V_jB_\ell\bigr)\\
&\le C\bigl(2^{(j-\ell)/2}+2^{(j-\ell)/2}\bigr)B_jB_\ell
\le C2^{-(\ell-j)/2}B_jB_\ell.
\end{aligned}
\]

The last line uses $V_\ell\le C2^{-\ell/2}B_\ell$ and $V_j\le C2^{-j/2}B_j$. Summing over the gap $d=\ell-j\ge4$ yields

\begin{equation}
\begin{aligned}
\sum_{\ell\ge j+4}\norm{\mathcal R(r_jr_\ell)}_2
&\le C\sum_{d=4}^\infty2^{-d/2}\sum_jB_jB_{j+d}\\
&\le C\sum_{d=4}^\infty2^{-d/2}\sum_jB_j^2
\le C\norm{\nabla r}_2\norm s_2.
\end{aligned}
\tag{A.8}
\end{equation}

Here the second step is Cauchy--Schwarz, and the third follows from the geometric series and (A.6). Interchanging $j$ and $\ell$ gives the same estimate for $j\ge\ell+4$. Together with Step 3, these cases cover every pair of integers.

\step{Step 5: summation and the limit.}
By (A.7)--(A.8) and symmetry,

\begin{equation}
\sum_{j,\ell\in\mathbb Z}\norm{\mathcal R(r_jr_\ell)}_2
\le C\norm{\nabla r}_2\norm s_2.
\tag{A.9}
\end{equation}

Let $r^{(N)}=\sum_{j=-N}^Nr_j$. Expanding the finite sum gives

\[
\norm{\mathcal R((r^{(N)})^2)}_2
\le\sum_{|j|,|\ell|\le N}\norm{\mathcal R(r_jr_\ell)}_2
\le C\norm{\nabla r}_2\norm s_2.
\]

Step 1 gives $r^{(N)}\to r$ in $H^1(\RR^2)$. By the two-dimensional Sobolev embedding $H^1\hookrightarrow L^4$ (see \cite[Section~1.4.2]{BCD2011}),

\[
\norm{(r^{(N)})^2-r^2}_2
\le\norm{r^{(N)}-r}_4
     \bigl(\norm{r^{(N)}}_4+\norm r_4\bigr)\longrightarrow0.
\]

The $L^2$ boundedness of $\mathcal R$ then implies $\mathcal R((r^{(N)})^2)\to\mathcal R(r^2)$ in $L^2$. Passing to the limit proves (A.3) for the stated $H^1$ divergence-free fields.\qedbox

\result{Remark A.2 }
The same dyadic argument gives, for $(r,s)\in H^1(\RR^2;\RR^2)$ with $r_x+s_y=0$,

\[
\norm{R_x^2(r^2)}_2\le C\norm{\nabla r}_2\norm s_2.
\]

Interchanging both the coordinates and the components gives

\[
\norm{R_y^2(s^2)}_2\le C\norm{\nabla s}_2\norm r_2.
\]

Both constants are absolute.
\section{Weighted estimates for the pressure combination}
All fields, Riesz transforms, and unqualified norms in this appendix are on the whole plane. Set

\[
\chi_\rho(y)=\frac\rho{\rho+|y|},\qquad\vartheta_\rho(y)=\chi_\rho(y)^2,\qquad\rho>0.
\]

We first prove a weighted estimate for odd harmonic functions, then combine it with Appendix A to treat the tangential square. Finally, we establish weighted estimates for the Riesz transforms. For the Muckenhoupt $A_2$ condition and the classical weighted Calder\'on--Zygmund theory, see \cite{Grafakos2014}. Since $\vartheta_\rho(y)$ is not an $A_2$ weight, that theory does not apply directly.
Instead, we exploit cancellation arising from reflection to derive
kernel estimates and obtain the required weighted $L^2$ bounds.

\result{Lemma B.1 }
Let $H\in L^2(\RR^2)$ be odd in $y$ and harmonic in $|y|<h$. Then

\begin{equation}
\norm{\chi_\rho H}_2\le C\chi_\rho(h)\norm H_2.
\tag{B.1}
\end{equation}

The constant is independent of $\rho$, $h$, and $H$.

\step{Proof. Step 1: the Fourier modes in the strip.}
After horizontal Fourier transformation, harmonicity becomes

\[
\partial_y^2\widehat H(k,y)-k^2\widehat H(k,y)=0.
\]

For almost every $k\ne0$, the general solution is

\[
\widehat H(k,y)=\alpha(k)\cosh(|k|y)+\beta(k)\sinh(|k|y).
\]

Oddness eliminates the $\cosh$ term, so

\[
\widehat H(k,y)=\beta(k)\sinh(|k|y).
\]

Thus the boundary value is zero, and only the growing $\sinh$ mode remains inside the strip.

\step{Step 2: cubic decay of the mass near the boundary.}
For $t>0$,

\[
\frac{\mathrm d}{\mathrm dt}\frac{\sinh t}{t}
=\frac{t\cosh t-\sinh t}{t^2}\ge0,
\]
 Consequently, for $0<y<z<h$,

\[
|\widehat H(k,y)|^2\le\frac{y^2}{z^2}|\widehat H(k,z)|^2.
\]

Set $g(y)=\|H(\cdot,y)\|_{L_x^2}^2$. Plancherel gives $g(y)\le(y/z)^2g(z)$. Writing $I(t)=\int_0^tg(y)\dd y$, we obtain, for $0<t<h$,

\[
I(t)=t\int_0^1g(tz)\dd z
\le t\left(\frac th\right)^2\int_0^1g(hz)\dd z
=\left(\frac th\right)^3I(h).
\]

\step{Step 3: insertion of the weight.}
For $y>0$,

\[
-\vartheta_\rho'(y)=\frac{2\vartheta_\rho(y)}{\rho+y}
\le\frac{2\vartheta_\rho(y)}y.
\]

Integration by parts and Step 2 give

\[
\begin{aligned}
\int_0^h\vartheta_\rho(y)g(y)\dd y
&=\vartheta_\rho(h)I(h)+\int_0^h[-\vartheta_\rho'(y)]I(y)\dd y\\
&\le\vartheta_\rho(h)I(h)
+\frac{2I(h)}{h^3}\int_0^h\vartheta_\rho(y)y^2\dd y.
\end{aligned}
\]

Since $y/(\rho+y)$ is increasing, $\vartheta_\rho(y)y^2\le\vartheta_\rho(h)h^2$ for $y<h$. Therefore

\[
\int_0^h\vartheta_\rho g\le3\vartheta_\rho(h)\int_0^hg.
\]

On $y\ge h$, use $\vartheta_\rho(y)\le\vartheta_\rho(h)$ directly. It follows that

\[
\int_0^\infty\vartheta_\rho g\le3\vartheta_\rho(h)\int_0^\infty g.
\]

Evenness of $|H|^2$ allows us to double both sides and take square roots, proving (B.1). \qedbox

\result{Lemma B.2 }
Let $A=(r,s)$ be the whole-plane extension in Lemma 5.3, so that

\[
A\in H^1_{\rm loc}(\RR^2;\RR^2)\cap L^\infty(\RR^2;\RR^2),\qquad
\nabla A\in L^2(\RR^2),\qquad r_x+s_y=0,
\]

with zero trace on $y=0$, $r$ odd in $y$, and $s$ even in $y$. Then

\begin{equation}
\norm{\chi_\rho R_xR_y(r^2)}_2
\le C\norm{\nabla A}_2\norm{\chi_\rho s}_2.
\tag{B.2}
\end{equation}

\step{Proof.}
We decompose the stream function by height so that each localized field remains divergence-free. Lemma A.2 estimates each localized square, and Lemma B.1 supplies the weight for sources away from the boundary. We then sum and identify the limit.

\step{Step 1: the stream function and Hardy bounds.}
Let $\Psi(x,y)=\int_0^yr(x,t)\dd t$. Since $r$ is odd, $\Psi$ is even. Applying the Hardy inequalities of Lemma 2.1 separately in the two half-planes gives

\[
\norm{A/y}_2\le C\norm{\nabla A}_2,\qquad
\norm{\Psi/y^2}_2\le\frac23\norm{r/y}_2.
\]

In particular, $\|\chi_\rho s\|_2\le\rho\|s/y\|_2<\infty$.

\step{Step 2: localization of the stream function.}
Set $T_j=2^j\rho$. Choose a real, even, smooth partition of squares $\sum_{j\ge0}\theta_j^2=1$ such that

\[
\begin{gathered}
\operatorname{supp}\theta_0\subset\{|y|\le2\rho\},\qquad
\operatorname{supp}\theta_j\subset\{T_j/2\le|y|\le2T_j\}\quad(j\ge1),\\
|\theta_j^{(\ell)}|\le C_\ell T_j^{-\ell}\quad(\ell=1,2).
\end{gathered}
\]

The supports have uniformly finite overlap. Denote them by $D_j$ and define

\[
\begin{aligned}
r_j&=\partial_y(\theta_j\Psi)=\theta_jr+\theta_j'\Psi,\\
s_j&=-\partial_x(\theta_j\Psi)=\theta_js,\qquad q_j=\theta_j'\Psi.
\end{aligned}
\]

 Write

\[
\begin{aligned}
E_j&=\norm{\nabla A}_{L^2(D_j)}
+T_j^{-1}\norm A_{L^2(D_j)}+T_j^{-2}\norm\Psi_{L^2(D_j)},\\
V_j&=\norm s_{L^2(D_j)}.
\end{aligned}
\]

On $D_j$, the cutoff scales satisfy $T_j^{-1}\le C/|y|$ and $T_j^{-2}\le C/y^2$. Step 1 and finite overlap imply

\[
\sum_jE_j^2\le C\bigl(\norm{\nabla A}_2^2+\norm{A/y}_2^2+\norm{\Psi/y^2}_2^2\bigr)
\le C\norm{\nabla A}_2^2.
\]

Since $\chi_\rho(y)\asymp \chi_\rho(T_j)$ on $D_j$,

\begin{equation}
\sum_jE_j^2\le C\norm{\nabla A}_2^2,\qquad
\sum_j\chi_\rho(T_j)^2V_j^2\le C\norm{\chi_\rho s}_2^2.
\tag{B.3}
\end{equation}

\step{Step 3: the localized squares and correction terms.}
Differentiating term by term,

\[
\begin{aligned}
\partial_xr_j&=\theta_jr_x-\theta_j's,\\
\partial_yr_j&=\theta_jr_y+2\theta_j'r+\theta_j''\Psi,\\
\partial_xq_j&=-\theta_j's.
\end{aligned}
\]

and hence

\[
\norm{\nabla r_j}_2\le CE_j,\quad \norm{s_j}_2\le V_j,\quad
\norm{q_j}_2\le CT_jE_j,\quad
\norm{\partial_xq_j}_2\le CT_j^{-1}V_j.
\]

Finite-height support and the Hardy bounds imply $(r_j,s_j)\in H^1(\RR^2;\RR^2)$. Lemma A.2 gives

\[
\norm{R_xR_y(r_j^2)}_2
\le C\norm{\nabla r_j}_2\norm{s_j}_2\le CE_jV_j.
\]

For the correction $q_j^2$, Lemma A.1 gives

\[
\norm{R_xR_y(q_j^2)}_2
\le C\norm{q_j}_2\norm{\partial_xq_j}_2\le CE_jV_j.
\]

\step{Step 4: the weight from odd harmonicity.}
Set $H_j=R_xR_y(r_j^2-q_j^2)$. Since $\theta_j$ and $\Psi$ are even and $\theta_j'$ is odd,

\[
r_j,q_j\text{ are odd},\qquad s_j\text{ is even}.
\]

Thus $r_j^2-q_j^2$ is even, and $H_j$ is odd. For $j\ge1$, the source vanishes on $|y|<T_j/2$. Since

\[
-\Delta H_j=\partial_x\partial_y(r_j^2-q_j^2),
\]

$H_j$ is harmonic in that strip. Lemma B.1 therefore yields

\begin{equation}
\norm{\chi_\rho H_j}_2
\le C\chi_\rho(T_j/2)\norm{H_j}_2
\le C\chi_\rho(T_j)E_jV_j.
\tag{B.4}
\end{equation}

For $j=0$, the same bound follows from $\chi_\rho\le1=2\chi_\rho(T_0)$. This is where parity is used: the odd harmonic output supplies the factor $\chi_\rho(T_j)$.

\step{Step 5: summation over height.}
Cauchy--Schwarz and (B.3) give

\[
\begin{aligned}
\sum_j\norm{\chi_\rho H_j}_2
&\le C\sum_j\chi_\rho(T_j)E_jV_j\\
&\le C\left(\sum_jE_j^2\right)^{1/2}
       \left(\sum_j\chi_\rho(T_j)^2V_j^2\right)^{1/2}\\
&\le C\norm{\nabla A}_2\norm{\chi_\rho s}_2.
\end{aligned}
\]

\step{Step 6: identification of the sum.}
By definition,

\[
r_j^2-q_j^2=\theta_j^2r^2+2\theta_j\theta_j'r\Psi.
\]

Since $\sum_j\theta_j^2=1$ and $2\sum_j\theta_j\theta_j'=\partial_y\sum_j\theta_j^2=0$, we have pointwise

\[
\sum_j(r_j^2-q_j^2)=r^2.
\]

To identify the limit after applying the Riesz transform, write

\[
\sigma_N=\sum_{j=0}^N\theta_j^2,\qquad
S_N=\sum_{j=0}^N(r_j^2-q_j^2)=\sigma_Nr^2+\sigma_N'r\Psi.
\]

Here $\sigma_N=1$ on $|y|<T_N$ and vanishes for $|y|>2T_N$, with $|\sigma_N'|\le CT_N^{-1}$ and $|\sigma_N''|\le CT_N^{-2}$. Expanding the gradient,

\[
\begin{aligned}
\nabla(S_N-r^2)={}&2(\sigma_N-1)r\nabla r+\sigma_N'r^2e_y\\
&+\sigma_N''r\Psi e_y+\sigma_N'\Psi\nabla r+\sigma_N'r\nabla\Psi.
\end{aligned}
\]

Using $|\Psi(x,y)|\le|y|\|A\|_\infty$ and $|\nabla\Psi|=|A|$, we estimate each term by

\[
\|\nabla(S_N-r^2)\|_2
\le C\|A\|_\infty\bigl(\|\nabla A\|_{2,\{|y|\ge T_N\}}+
\|A/y\|_{2,\{|y|\ge T_N\}}\bigr)\longrightarrow0.
\]

Hence $\nabla R_xR_yS_N\to R_xR_y\nabla(r^2)$ in $L^2$. On the other hand, Step 5 gives a weighted $L^2$ limit of $R_xR_yS_N=\sum_{j=0}^NH_j$. This limit and $R_xR_y(r^2)$ have the same gradient and differ only by a constant. Both representatives are odd, so the constant is zero. Step 5 now proves (B.2).\qedbox

\pdfbookmark[2]{Lemma B.3: weighted Riesz estimates for odd outputs}{lemmaB3}
\result{Lemma B.3 }
Suppose $F$ is even and $H$ is odd in $y$, with $\chi_\rho F,\chi_\rho H\in L^2(\RR^2)$. Then

\begin{equation}
\begin{aligned}
\|\chi_\rho R_xR_yF\|_2&\le C\|\chi_\rho F\|_2,\\
\|\chi_\rho R_x^2H\|_2+\|\chi_\rho R_y^2H\|_2&\le C\|\chi_\rho H\|_2,
\end{aligned}\tag{B.5}
\end{equation}

where $C$ is independent of $\rho$.
\step{Proof. Step 1: the reflected kernel.}
We first estimate $R_xR_y$ on even inputs and $R_y^2-R_x^2$ on odd inputs, and then separate the two pure transforms by a multiplier identity. Let $\eta>0$ be the output height, $y>0$ the source height, and $s$ the horizontal source coordinate. Set $T_1=R_xR_y$ and $T_2=R_y^2-R_x^2$. These are second-order Riesz transforms; see \cite{Stein1970} for their singular-integral representation. Away from the origin, their kernels are

\[
K_1(\xi,t)=\frac{\xi t}{\pi(\xi^2+t^2)^2},\qquad
K_2(\xi,t)=\frac{t^2-\xi^2}{\pi(\xi^2+t^2)^2}.
\]

Begin with ordinary $L^2$ inputs. Split the whole-plane integral into the two half-planes, reflect the lower half, and use the parity of the input. This gives

\[
\begin{aligned}
T_1F(x,\eta)&=\operatorname{p.v.}\int_0^\infty\int_\RR
\mathcal K_1(x-s;\eta,y)F(s,y)\,ds\,dy,\\
T_2H(x,\eta)&=\operatorname{p.v.}\int_0^\infty\int_\RR
\mathcal K_2(x-s;\eta,y)H(s,y)\,ds\,dy,
\end{aligned}
\]

where

\[
\begin{aligned}
\mathcal K_1(\xi;\eta,y)&=K_1(\xi,\eta-y)+K_1(\xi,\eta+y),\\
\mathcal K_2(\xi;\eta,y)&=K_2(\xi,\eta-y)-K_2(\xi,\eta+y).
\end{aligned}
\]

 Both outputs are odd in the normal variable.

\step{Step 2: estimates for separated heights.}
The kernel formulas and the substitution $\xi=|t|z$ give, for $i=1,2$ and $t\ne0$,

\[
\|K_i(\cdot,t)\|_{L^1_\xi}\le C|t|^{-1},\qquad
\|\partial_tK_i(\cdot,t)\|_{L^1_\xi}\le C|t|^{-2}.
\]

If $y\ge4\eta$, oddness of $K_1$ and evenness of $K_2$ in $t$ turn the paired kernels into

\[
\begin{aligned}
\mathcal K_1(\xi;\eta,y)&=K_1(\xi,y+\eta)-K_1(\xi,y-\eta),\\
\mathcal K_2(\xi;\eta,y)&=K_2(\xi,y-\eta)-K_2(\xi,y+\eta).
\end{aligned}
\]

The fundamental theorem of calculus, with $t\asymp y$ on $[y-\eta,y+\eta]$, gives

\[
\|\mathcal K_i(\cdot;\eta,y)\|_{L^1_\xi}
\le\int_{y-\eta}^{y+\eta}\|\partial_tK_i(\cdot,t)\|_{L^1_\xi}\,dt
\le C\frac\eta{y^2}.
\]

If instead $\eta\ge4y$, then $\eta-y\asymp\eta+y\asymp\eta$, and the ordinary kernel bound gives

\[
\|\mathcal K_i(\cdot;\eta,y)\|_{L^1_\xi}
\le\|K_i(\cdot,\eta-y)\|_1+\|K_i(\cdot,\eta+y)\|_1
\le\frac C\eta.
\]

In the first case, reflection improves $1/y$ to $\eta/y^2$. This extra factor makes the weighted interaction between distant layers summable.

\step{Step 3: height layers and neighboring interactions.}
Let $T_j=2^j\rho$ and define

\[
J_0=(0,2\rho),\qquad J_j=[2^j\rho,2^{j+1}\rho)\quad(j\ge1).
\]

Then $|J_j|\asymp T_j$, and $\chi_\rho(y)\asymp \chi_\rho(T_j)$ on $J_j$. In what follows, $T$ denotes $T_1$ or $T_2$, and $f$ has the corresponding even or odd parity. Set $f_j=1_{J_j}f$ on the upper half-plane and extend it with the same parity as $f$. Write

\[
d_j=\|f_j\|_{L^2(\RR\times J_j)},\qquad
b_j=\chi_\rho(T_j)d_j.
\]

First consider finitely many source layers. For $|j-k|\le3$, unweighted Riesz-transform boundedness gives

\[
\|Tf_j\|_{L^2(\RR\times J_k)}\le C\|f_j\|_{L^2(\RR^2)}
\le Cd_j.
\]

The weights on neighboring layers are comparable, so

\[
\chi_\rho(T_k)\|Tf_j\|_{L^2(\RR\times J_k)}\le Cb_j.
\]

No absolute integration of the singular kernel is used here; the singularity on the same or neighboring layers is handled by unweighted $L^2$ boundedness.

\step{Step 4: the source is higher, $j\ge k+4$.}
For $y\in J_j$ and $\eta\in J_k$, we have $y\ge4\eta$. The paired-kernel representation and Young's inequality in the horizontal variable give, at fixed $\eta$,

\[
\|Tf_j(\cdot,\eta)\|_{L^2_x}
\le C\frac{T_k}{T_j^2}\int_{J_j}\|f_j(\cdot,y)\|_{L^2_x}\,dy
\le C\frac{T_k}{T_j^2}|J_j|^{1/2}d_j.
\]

Taking the $L^2$ norm over $\eta\in J_k$ yields

\[
\|Tf_j\|_{L^2(\RR\times J_k)}
\le C\frac{T_k}{T_j^2}|J_j|^{1/2}|J_k|^{1/2}d_j
\le C\left(\frac{T_k}{T_j}\right)^{3/2}d_j.
\]

Since $\chi_\rho(T_k)/\chi_\rho(T_j)\le CT_j/T_k$,

\[
\chi_\rho(T_k)\|Tf_j\|_{L^2(\RR\times J_k)}
\le C2^{-(j-k)/2}b_j.
\]

\step{Step 5: the output is higher, $k\ge j+4$.}
Now $\eta\ge4y$. Apply horizontal Young, then Cauchy--Schwarz in the source height, and finally integrate over the output height to obtain

\[
\begin{aligned}
\|Tf_j(\cdot,\eta)\|_{L^2_x}&\le CT_k^{-1}|J_j|^{1/2}d_j,\\
\|Tf_j\|_{L^2(\RR\times J_k)}&\le CT_k^{-1}|J_j|^{1/2}|J_k|^{1/2}d_j
\le C\left(\frac{T_j}{T_k}\right)^{1/2}d_j.
\end{aligned}
\]

Using $\chi_\rho(T_k)/\chi_\rho(T_j)\le CT_j/T_k$ gives

\[
\chi_\rho(T_k)\|Tf_j\|_{L^2(\RR\times J_k)}\le C2^{-3(k-j)/2}b_j.
\]

\step{Step 6: summation and general weighted inputs.}
Steps 3--5 imply

\[
\chi_\rho(T_k)\|Tf\|_{L^2(\RR\times J_k)}\le C\sum_{j\ge0}c_{k-j}b_j,
\qquad
c_m=\begin{cases}
1,&|m|\le3,\\
2^{m/2},&m\le-4,\\
2^{-3m/2},&m\ge4.
\end{cases}
\]

Since $\sum_{m\in\mathbb Z}c_m<\infty$, Young's convolution inequality on $\mathbb Z$ (see \cite[Section~1.2]{Grafakos2014}) and comparability of the weights within each layer give

\[
\|\chi_\rho Tf\|_{L^2(\eta>0)}^2
\le C\sum_{k\ge0}\left(\sum_{j\ge0}c_{k-j}b_j\right)^2
\le C\sum_{j\ge0}b_j^2
\le C\|\chi_\rho f\|_{L^2(y>0)}^2.
\]

The squared moduli of both input and output are even, so this proves the whole-plane estimate. For a general input, let $f^{(N)}=\sum_{j=0}^Nf_j$. Then $f^{(N)}\in L^2$ and $\chi_\rho f^{(N)}\to \chi_\rho f$ in $L^2$. Applying the same estimate to differences makes $\chi_\rho Tf^{(N)}$ Cauchy. Its limit defines the weighted operator independently of the approximation. We have proved the first estimate in (B.5), as well as $\|\chi_\rho(R_y^2-R_x^2)H\|_2\le C\|\chi_\rho H\|_2$.

\step{Step 7: the pure second-order transforms.}
The identity $R_x^2+R_y^2=-I$ gives

\[
R_x^2=-\tfrac12 I-\tfrac12(R_y^2-R_x^2),\qquad
R_y^2=-\tfrac12 I+\tfrac12(R_y^2-R_x^2).
\]

and hence

\[
\|\chi_\rho R_x^2H\|_2+\|\chi_\rho R_y^2H\|_2
\le\|\chi_\rho H\|_2+\|\chi_\rho(R_y^2-R_x^2)H\|_2
\le C\|\chi_\rho H\|_2.
\]

\qedbox

\step{\textbf{Proof of Lemma 5.3.}}
Use the whole-plane extensions and norm convention above. Since $r$ is odd and $s$ is even,

\[
r^2,s^2\text{ are even},\qquad rs\text{ is odd}.
\]

By definition,

\[
G(A\otimes A)=R_xR_y(r^2)-R_xR_y(s^2)+(R_y^2-R_x^2)(rs).
\]

The term without a factor of $s$ is controlled by Lemma B.2:

\[
\norm{\chi_\rho R_xR_y(r^2)}_2
\le C\norm{\nabla A}_2\norm{\chi_\rho s}_2.
\]

Since $s^2$ is even, the first estimate in Lemma B.3 gives

\[
\begin{aligned}
\norm{\chi_\rho R_xR_y(s^2)}_2
&\le C\norm{\chi_\rho s^2}_2\\
&\le C\norm s_\infty\norm{\chi_\rho s}_2
\le C\norm A_\infty\norm{\chi_\rho s}_2.
\end{aligned}
\]

Finally, $rs$ is odd, so the second estimate in Lemma B.3 gives

\[
\begin{aligned}
\norm{\chi_\rho(R_y^2-R_x^2)(rs)}_2
&\le C\norm{\chi_\rho rs}_2\\
&\le C\norm r_\infty\norm{\chi_\rho s}_2
\le C\norm A_\infty\norm{\chi_\rho s}_2.
\end{aligned}
\]

Adding these three estimates yields

\[
\begin{aligned}
\norm{\chi_\rho G(A\otimes A)}_2
&\le\norm{\chi_\rho R_xR_y(r^2)}_2+\norm{\chi_\rho R_xR_y(s^2)}_2\\
&\quad+\norm{\chi_\rho(R_y^2-R_x^2)(rs)}_2\\
&\le C\bigl(\norm{\nabla A}_2+\norm A_\infty\bigr)\norm{\chi_\rho s}_2.
\end{aligned}
\]

\qedbox

\section*{Acknowledgements}

During the preparation of this manuscript, the author used
ChatGPT 6 Astra to assist with literature searches, the review
of mathematical arguments, and language editing.
The author has checked and substantially revised all AI-generated
mathematical arguments and takes full responsibility for the
content of the final manuscript.


\begin{thebibliography}{99}
\bibitem{Amick1988}
C. J. Amick, On Leray's problem of steady Navier--Stokes flow past a body in the plane, \emph{Acta Math.} \textbf{161} (1988), 71--130.

\bibitem{BCD2011}
H. Bahouri, J.-Y. Chemin, and R. Danchin, \emph{Fourier Analysis and Nonlinear Partial Differential Equations}, Grundlehren der mathematischen Wissenschaften, vol.~343, Springer, Heidelberg, 2011. \href{https://doi.org/10.1007/978-3-642-16830-7}{doi:10.1007/978-3-642-16830-7}.

\bibitem{BFZ2013}
M. Bildhauer, M. Fuchs, and G. Zhang, Liouville-type theorems for steady flows of degenerate power law fluids in the plane, \emph{J. Math. Fluid Mech.} \textbf{15} (2013), no. 3, 583--616. \href{https://doi.org/10.1007/s00021-012-0122-1}{doi:10.1007/s00021-012-0122-1}.

\bibitem{Bony1981}
J.-M. Bony, Calcul symbolique et propagation des singularit\'es pour les \'equations aux d\'eriv\'ees partielles non lin\'eaires, \emph{Ann. Sci. \'Ecole Norm. Sup. (4)} \textbf{14} (1981), no.~2, 209--246. \href{https://doi.org/10.24033/asens.1404}{doi:10.24033/asens.1404}.

\bibitem{CLP2024}
A. D. Calderon, V. Le, and T. Phan, On stationary Navier--Stokes equations in the upper-half plane, \emph{Acta Appl. Math.} \textbf{189} (2024), Paper No. 7. \href{https://doi.org/10.1007/s10440-024-00636-3}{doi:10.1007/s10440-024-00636-3}.

\bibitem{Fuchs2012}
M. Fuchs, Liouville theorems for stationary flows of shear thickening fluids in the plane, \emph{J. Math. Fluid Mech.} \textbf{14} (2012), no. 3, 421--444. \href{https://doi.org/10.1007/s00021-011-0070-1}{doi:10.1007/s00021-011-0070-1}.

\bibitem{FuchsZhang2012}
M. Fuchs and G. Zhang, Liouville theorems for entire local minimizers of energies defined on the class $L\log L$ and for entire solutions of the stationary Prandtl--Eyring fluid model, \emph{Calc. Var. Partial Differential Equations} \textbf{44} (2012), nos. 1--2, 271--295. \href{https://doi.org/10.1007/s00526-011-0434-7}{doi:10.1007/s00526-011-0434-7}.

\bibitem{FZ2011}
M. Fuchs and X. Zhong, A note on a Liouville type result of Gilbarg and Weinberger for the stationary Navier--Stokes equations in 2D, \emph{J. Math. Sci. (N.Y.)} \textbf{178} (2011), no. 6, 695--703. \href{https://doi.org/10.1007/s10958-011-0578-1}{doi:10.1007/s10958-011-0578-1}.

\bibitem{Galdi2011}
G. P. Galdi, \emph{An Introduction to the Mathematical Theory of the Navier--Stokes Equations: Steady-State Problems}, second ed., Springer Monographs in Mathematics, Springer, New York, 2011. \href{https://doi.org/10.1007/978-0-387-09620-9}{doi:10.1007/978-0-387-09620-9}.

\bibitem{GHM2014}
Y. Giga, P.-Y. Hsu, and Y. Maekawa, A Liouville theorem for the planer Navier--Stokes equations with the no-slip boundary condition and its application to a geometric regularity criterion, \emph{Comm. Partial Differential Equations} \textbf{39} (2014), no. 10, 1906--1935. \href{https://doi.org/10.1080/03605302.2014.912662}{doi:10.1080/03605302.2014.912662}.

\bibitem{GW1978}
D. Gilbarg and H. F. Weinberger, Asymptotic properties of steady plane solutions of the Navier--Stokes equations with bounded Dirichlet integral, \emph{Ann. Scuola Norm. Sup. Pisa Cl. Sci. (4)} \textbf{5} (1978), no. 2, 381--404. \href{https://www.numdam.org/item/ASNSP_1978_4_5_2_381_0/}{Numdam}.

\bibitem{Grafakos2014}
L. Grafakos, \emph{Classical Fourier Analysis}, third ed., Graduate Texts in Mathematics, vol.~249, Springer, New York, 2014. \href{https://doi.org/10.1007/978-1-4939-1194-3}{doi:10.1007/978-1-4939-1194-3}.

\bibitem{KNSS2009}
G. Koch, N. Nadirashvili, G. A. Seregin, and V. \v Sver\'ak, Liouville theorems for the Navier--Stokes equations and applications, \emph{Acta Math.} \textbf{203} (2009), no. 1, 83--105. \href{https://arxiv.org/abs/0709.3599}{arXiv:0709.3599}.

\bibitem{KPR2019}
M. V. Korobkov, K. Pileckas, and R. Russo, On convergence of arbitrary $D$-solution of steady Navier--Stokes system in 2D exterior domains, \emph{Arch. Ration. Mech. Anal.} \textbf{233} (2019), 385--407. \href{https://doi.org/10.1007/s00205-019-01359-8}{doi:10.1007/s00205-019-01359-8}.

\bibitem{KR2022}
M. Korobkov and X. Ren, Leray's plane stationary solutions at small Reynolds numbers, \emph{J. Math. Pures Appl. (9)} \textbf{158} (2022), 71--89. \href{https://doi.org/10.1016/j.matpur.2021.08.010}{doi:10.1016/j.matpur.2021.08.010}.

\bibitem{KR2025}
M. Korobkov and X. Ren, On basic velocity estimates for the plane steady-state Navier--Stokes system and its applications, \emph{J. Math. Fluid Mech.} \textbf{27} (2025), Paper No. 37. \href{https://doi.org/10.1007/s00021-025-00939-x}{doi:10.1007/s00021-025-00939-x}. Preprint: \href{https://arxiv.org/abs/2405.17884v1}{arXiv:2405.17884v1}; the numbering of Corollary 5 and equation (1.9) cited here refers to this version.

\bibitem{Leray1933}
J. Leray, \'Etude de diverses \'equations int\'egrales non lin\'eaires et de quelques probl\`emes que pose l'hydrodynamique, \emph{J. Math. Pures Appl. (9)} \textbf{12} (1933), 1--82. \href{https://www.numdam.org/item/JMPA_1933_9_12__1_0/}{Numdam}.

\bibitem{Leray1934}
J. Leray, Essai sur les mouvements plans d'un liquide visqueux que limitent des parois, \emph{J. Math. Pures Appl. (9)} \textbf{13} (1934), 331--418. \href{https://www.numdam.org/item/JMPA_1934_9_13__331_0/}{Numdam}.

\bibitem{Sazonov1999}
L. I. Sazonov, On the asymptotic behavior of the solution of the two-dimensional stationary problem of the flow past a body far from it, \emph{Math. Notes} \textbf{65} (1999), no. 2, 202--207. \href{https://doi.org/10.1007/BF02679817}{doi:10.1007/BF02679817}.

\bibitem{Seregin2014}
G. Seregin, Liouville theorem for 2D Navier--Stokes equations in half space, \emph{Zap. Nauchn. Sem. POMI} \textbf{425} (2014), 137--148; English translation in \emph{J. Math. Sci. (N.Y.)} \textbf{210} (2015), no. 6, 849--856. \href{https://doi.org/10.1007/s10958-015-2595-y}{doi:10.1007/s10958-015-2595-y}.

\bibitem{Smith1965}
D. R. Smith, Estimates at infinity for stationary solutions of the Navier--Stokes equations in two dimensions, \emph{Arch. Ration. Mech. Anal.} \textbf{20} (1965), 341--372. \href{https://doi.org/10.1007/BF00282357}{doi:10.1007/BF00282357}.

\bibitem{Stein1970}
E. M. Stein, \emph{Singular Integrals and Differentiability Properties of Functions}, Princeton Mathematical Series, vol.~30, Princeton University Press, Princeton, NJ, 1970. \href{https://www.jstor.org/stable/j.ctt1bpmb07}{Publisher's edition on JSTOR}.

\bibitem{WW2022}
L. Wang and W. Wang, Asymptotic properties of steady plane solutions of the Navier--Stokes equations in the exterior of a half-space, preprint (2022), \href{https://arxiv.org/abs/2208.01458v2}{arXiv:2208.01458v2}.

\bibitem{WW2023}
L. Wang and W. Wang, Asymptotic properties of steady plane solutions of the Navier--Stokes equations in a cone-like domain, \emph{J. Math. Fluid Mech.} \textbf{25} (2023), Paper No. 75. \href{https://doi.org/10.1007/s00021-023-00818-3}{doi:10.1007/s00021-023-00818-3}.

\bibitem{Zhang2015}
G. Zhang, Liouville theorems for stationary flows of shear thickening fluids in 2D, \emph{Ann. Acad. Sci. Fenn. Math.} \textbf{40} (2015), no. 2, 889--905. \href{https://doi.org/10.5186/aasfm.2015.4052}{doi:10.5186/aasfm.2015.4052}.

\bibitem{BW2012}
C.~Boeckle and P.~Wittwer,
Decay estimates for steady solutions of the Navier--Stokes equations
in two dimensions in the presence of a wall,
\textit{SIAM J. Math. Anal.} \textbf{44} (2012), no.~5, 3346--3368.

\end{thebibliography}
\end{document}